\documentclass[titlepage,11pt]{article}

\def\doublespacing{\parskip 5 pt plus 1 pt \baselineskip 16pt
	\lineskip 13 pt \normallineskip 13 pt}
\usepackage{epsfig}
\usepackage{color,graphics}
\usepackage{pdfpages}
\usepackage{pictex}
\begingroup\makeatletter\ifx\SetFigFontNFSS\undefined%
\gdef\SetFigFontNFSS#1#2#3#4#5{%
\reset@font\fontsize{#1}{#2pt}%
\fontfamily{#3}\fontseries{#4}\fontshape{#5}%
\selectfont}%
\fi\endgroup%
\usepackage{amsmath}
\usepackage{amssymb,dsfont}
\usepackage{textcomp}
\usepackage{verbatim} 
\usepackage{color}

\newtheorem{thm}{Theorem}[section]
\newtheorem{prop}[thm]{Proposition}
\newtheorem{cor}[thm]{Corollary}
\newtheorem{lem}[thm]{Lemma}
\newtheorem{conj}[thm]{Conjecture}
\newtheorem{prob}[thm]{Problem}
\newtheorem{exa}[thm]{Example}
\newtheorem{obs}[thm]{Observation}
\newtheorem{defn}[thm]{Definition}
\newtheorem{alg}[thm]{Algorithm}
\newtheorem{que}[thm]{Question}

\newcommand{\ben}{\begin{enumerate}}
\newcommand{\een}{\end{enumerate}}
\newcommand{\bobs}{\begin{obs}}
\newcommand{\eobs}{\end{obs}}
\newcommand{\ble}{\begin{lem}}
\newcommand{\ele}{\end{lem}}
\newcommand{\bthm}{\begin{thm}}
\newcommand{\ethm}{\end{thm}}
\newcommand{\bpr}{\begin{prop}}
\newcommand{\epr}{\end{prop}}
\newcommand{\bco}{\begin{cor}}
\newcommand{\eco}{\end{cor}}
\newcommand{\bcon}{\begin{conj}}
\newcommand{\econ}{\end{conj}}
\newcommand{\bprob}{\begin{prob}}
\newcommand{\eprob}{\end{prob}}
\newcommand{\bde}{\begin{defn}}
\newcommand{\ede}{\end{defn}}
\newcommand{\bex}{\begin{exa}}
\newcommand{\eex}{\end{exa}}
\newcommand{\balg}{\begin{alg}}
\newcommand{\ealg}{\end{alg}}
\newcommand{\bque}{\begin{que}}
\newcommand{\eque}{\end{que}}

\newcommand{\barr}{\begin{array}}
\newcommand{\earr}{\end{array}}
\newcommand{\btab}{\begin{tabular}}
\newcommand{\etab}{\end{tabular}}
\newcommand{\beq}{\begin{equation}}
\newcommand{\eeq}{\end{equation}}
\newcommand{\bea}{\begin{eqnarray*}}
\newcommand{\eea}{\end{eqnarray*}}
\newcommand{\bce}{\begin{center}}
\newcommand{\ece}{\end{center}}
\newcommand{\bpi}{\begin{picture}}
\newcommand{\epi}{\end{picture}}
\newcommand{\bfi}{\begin{figure}[h] \begin{center}}
\newcommand{\efi}{\end{center} \end{figure}}
\newcommand{\capt}{\caption}
\newcommand{\bsl}{\begin{slide}{}}

\newcommand{\pf}{{\bf Proof.}}

\newcommand{\qed}{\rule{1ex}{1ex}}

\newcommand{\ul}{\underline}

\newcommand{\hs}[1]{\hspace{#1}}
\newcommand{\vs}[1]{\vspace{#1}}

\newcommand{\flf}[2]{\left\lfloor\frac{#1}{#2}\right\rfloor}
\newcommand{\cef}[2]{\left\lceil\frac{#1}{#2}\right\rceil}

\newcommand{\llra}{\longleftrightarrow}

\newcommand{\LRA}{\Leftrightarrow}

\newcommand{\al}{\alpha}
\newcommand{\be}{\beta}

\newcommand{\NN}{\mathbb{N}}
\newcommand{\ZZ}{\mathbb{Z}}

\newcommand{\diam}{\mathop{\rm diam}}

\newcommand{\ID}{\mathop{\rm ID}}

\newcommand \vd {\vec{d}}

\newcommand{\flo}[1]{\left\lfloor {#1} \right\rfloor}
\newcommand{\cei}[1]{\left\lceil {#1} \right\rceil}

\newcommand{\dist}[2]{d(u_{#1}, u_{#2})}
\newcommand{\dok}[1]{\dist 1{k_{#1}} }
\newcommand{\dkk}[1]{\dist k{k_{#1}} }
\newcommand{\dokp}[1]{\dist 1{k_{#1}'} }
\newcommand{\dkkp}[1]{\dist k{k_{#1}'} }

\newcommand{\ujo}{u_{j_1} }
\newcommand{\ujt}{u_{j_2} }
\newcommand{\flnt}  {\flo{\frac n2} }

\newcommand{\cent}  {\cei{\frac n2} }
\newcommand{\uel}{u_\ell}
\newcommand{\uelp}{u_{\ell'}}

\begin{document}

\bce
{\large \bf  A Note on ID-Colorings and Symmetric Colorings of Cycles}

\bigskip
Yuya Kono 

\ece
\vs{-1mm}

\bce
{\bf Abstract}
\ece
\vs{-5mm}
\begin{quote}
A red-white coloring of a nontrivial connected graph $G$ is an assignment of red and white colors to the vertices of~$G$.
Associated with each vertex $v$ of $G$ of diameter $d$ is a $d$-vector, called the code of $v$, whose $i$th coordinate is the number of red vertices at distance $i$ from $v$. 
A red-white coloring of $G$ for which distinct vertices have distinct codes is called an ID-coloring of $G$.
In 2025, a criterion to determine whether a red-white coloring of a path is an ID-coloring or not was presented by Kono, with the aid of a result shown by Marcelo et al. in 2024.
The criterion utilizes the fact that ID-colorings of paths are ``opposite'' of colorings with a certain symmetry.
In this paper, we establish a similar criterion that can be applied for cycles whose order is a prime number at least 3.
In order to do so, we employ an analogous approaches used for the criterion for paths, i.e., we pay attention to symmetries of given red-white colorings of cycles.
\end{quote} 

\medskip

\vs{.2cm}


\vs{1mm}

\doublespacing

\section{Introduction}

Let $G$ be a  connected graph of diameter $d\ge 2$ and let there be given a red-white vertex coloring $c$ of the graph~$G$ 
where  at least  one vertex  is colored red. That is,   the {\bf  color}~$c(v)$ of a vertex $v$ in~$G$ is either red or white and $c(v)$ is red for at least  one vertex~$v$ of~$G$. With each vertex $v$ of~$G$,  there is an associated {\bf $d$-vector $\vec{d}(v)=(a_1, a_2, \ldots, a_d)$} called the {\bf code}  of $v$ corresponding to~$c$,  where  the $i$th coordinate~$a_i$ is the number of red vertices at distance~$i$ from $v$ for $1 \le i \le d$. 
If  distinct vertices of~$G$ have distinct codes,  then~$c$ is called an 
{\bf identification coloring} or {\bf ID-coloring}.  Equivalently, 
an identification coloring of a connected graph~$G$  is an assignment of the color red to a nonempty subset of $V(G)$ (with the color white  assigned to the remaining vertices of~$G$) such that for every two vertices~$u$ and~$v$ of~$G$, there is an integer~$k$ with~$1 \le k \le d$ such that
the number of red vertices at distance~$k$ from~$u$ is different from 
the number of red vertices at distance~$k$ from~$v$.  A graph  possessing  an identification coloring  is an {\bf ID-graph}.
It is known that not all connected graphs are ID-graphs. 
The minimum number of red vertices among all ID-colorings of an ID-graph~$G$ is the 
{\bf identification number} or {\bf  ID-number} $\ID(G)$  of~$G$. 
This concept was introduced by Gary Chartrand and first studied in \cite{IDgraph1}.


In 2025, a criterion to determine whether a red-white coloring of a path is an ID-coloring or not was presented by Kono in \cite{KonoPath}, with the aid of a result shown by Marcelo et al. in 2024 in \cite{Marcelo}, where they focused on a symmetry of a given red-white coloring.
Let $P_n = (u_1, u_2, \cdots, u_n)$ be the path of order $n\ge 2$. 
We say that vertices $u_i$ and $u_j$ on $P_n$ are {\bf partners} if $i+j=n+1$. 
If $n$ is odd, then the vertex $u_{\cent}$
is called the {\bf central vertex} of $P_n$.
The partner of the central vertex of an odd path is the central vertex itself.
A red-white coloring of $P_n$ is called {\bf symmetric}  if each pair of partners of $P_n$ have the same color.
Any symmetric colorings of an odd path can assign white or red to the central vertex.
Figure \ref{P6-P9-pw} shows examples of symmetric colorings of $P_6$ and $P_9$, and each pair of partners is indicated by an arc.

\vs{2mm}
\bfi
\scalebox{0.8}{
{\unitlength 0.1in%
\begin{picture}(62.4800,11.5400)(5.1200,-17.2400)%
%
\special{sh 1.000}%
\special{ia 976 624 54 54 0.0000000 6.2831853}%
\special{pn 13}%
\special{ar 976 624 54 54 0.0000000 6.2831853}%
%
\special{pn 20}%
\special{pa 922 624}%
\special{pa 598 624}%
\special{fp}%
%
\special{sh 0}%
\special{ia 1396 624 54 54 0.0000000 6.2831853}%
\special{pn 13}%
\special{ar 1396 624 54 54 0.0000000 6.2831853}%
%
\special{pn 20}%
\special{pa 1342 624}%
\special{pa 1018 624}%
\special{fp}%
%
\special{sh 0}%
\special{ia 1816 624 54 54 0.0000000 6.2831853}%
\special{pn 13}%
\special{ar 1816 624 54 54 0.0000000 6.2831853}%
%
\special{pn 20}%
\special{pa 1762 624}%
\special{pa 1438 624}%
\special{fp}%
%
\special{sh 1.000}%
\special{ia 2236 624 54 54 0.0000000 6.2831853}%
\special{pn 13}%
\special{ar 2236 624 54 54 0.0000000 6.2831853}%
%
\special{pn 20}%
\special{pa 2182 624}%
\special{pa 1858 624}%
\special{fp}%
%
\special{sh 1.000}%
\special{ia 2656 624 54 54 0.0000000 6.2831853}%
\special{pn 13}%
\special{ar 2656 624 54 54 0.0000000 6.2831853}%
%
\special{pn 20}%
\special{pa 2602 624}%
\special{pa 2278 624}%
\special{fp}%
%
\special{sh 1.000}%
\special{ia 566 624 54 54 0.0000000 6.2831853}%
\special{pn 13}%
\special{ar 566 624 54 54 0.0000000 6.2831853}%
%
\special{sh 1.000}%
\special{ia 3776 624 54 54 0.0000000 6.2831853}%
\special{pn 13}%
\special{ar 3776 624 54 54 0.0000000 6.2831853}%
%
\special{pn 20}%
\special{pa 3722 624}%
\special{pa 3398 624}%
\special{fp}%
%
\special{sh 1.000}%
\special{ia 4196 624 54 54 0.0000000 6.2831853}%
\special{pn 13}%
\special{ar 4196 624 54 54 0.0000000 6.2831853}%
%
\special{pn 20}%
\special{pa 4142 624}%
\special{pa 3818 624}%
\special{fp}%
%
\special{sh 0}%
\special{ia 4616 624 54 54 0.0000000 6.2831853}%
\special{pn 13}%
\special{ar 4616 624 54 54 0.0000000 6.2831853}%
%
\special{pn 20}%
\special{pa 4562 624}%
\special{pa 4238 624}%
\special{fp}%
%
\special{sh 1.000}%
\special{ia 5036 624 54 54 0.0000000 6.2831853}%
\special{pn 13}%
\special{ar 5036 624 54 54 0.0000000 6.2831853}%
%
\special{pn 20}%
\special{pa 4982 624}%
\special{pa 4658 624}%
\special{fp}%
%
\special{sh 0}%
\special{ia 5456 624 54 54 0.0000000 6.2831853}%
\special{pn 13}%
\special{ar 5456 624 54 54 0.0000000 6.2831853}%
%
\special{pn 20}%
\special{pa 5402 624}%
\special{pa 5078 624}%
\special{fp}%
%
\special{sh 0}%
\special{ia 3366 624 54 54 0.0000000 6.2831853}%
\special{pn 13}%
\special{ar 3366 624 54 54 0.0000000 6.2831853}%
%
\special{sh 1.000}%
\special{ia 5866 624 54 54 0.0000000 6.2831853}%
\special{pn 13}%
\special{ar 5866 624 54 54 0.0000000 6.2831853}%
%
\special{pn 20}%
\special{pa 5812 624}%
\special{pa 5488 624}%
\special{fp}%
%
\special{sh 1.000}%
\special{ia 6286 624 54 54 0.0000000 6.2831853}%
\special{pn 13}%
\special{ar 6286 624 54 54 0.0000000 6.2831853}%
%
\special{pn 20}%
\special{pa 6232 624}%
\special{pa 5908 624}%
\special{fp}%
%
\special{sh 0}%
\special{ia 6706 624 54 54 0.0000000 6.2831853}%
\special{pn 13}%
\special{ar 6706 624 54 54 0.0000000 6.2831853}%
%
\special{pn 20}%
\special{pa 6652 624}%
\special{pa 6328 624}%
\special{fp}%
%
\special{pn 13}%
\special{ar 5025 770 845 495 6.2831853 3.1415927}%
%
\special{pn 13}%
\special{ar 5030 770 1270 744 6.2831853 3.1415927}%
%
\special{pn 13}%
\special{ar 5030 740 1680 984 6.2831853 3.1415927}%
%
\special{pn 13}%
\special{ar 5015 740 425 212 6.2831853 3.1415927}%
%
\special{pn 13}%
\special{ar 1615 740 625 312 6.2831853 3.1415927}%
%
\special{pn 13}%
\special{ar 1610 740 1050 524 6.2831853 3.1415927}%
%
\special{pn 13}%
\special{ar 1585 720 215 172 6.2831853 3.1415927}%
\end{picture}}%
}
\capt{Symmetric colorings of $P_6$ and $P_9$}
\label{P6-P9-pw}
\efi

The following theorem was established in 2025. The sufficient condition was shown by Kono in \cite{KonoPath} and the necessary condition was proven by Marcelo et al. in \cite{Marcelo}.

\bthm
\label{equiv}
Let $n\ge 2$ and let $c$ be a red-white coloring of the path $P_n$ under which the end vertices of $P_n$ are colored red. The coloring $c$ is an ID-coloring if and only if $c$ is not symmetric.
\ethm

Utilizing this theorem, the following criterion for paths was presented in \cite{KonoPath}.

\bthm \label{restrictid}
Let $r \ge 2$ and $n \ge r$. Suppose that $c$ is a red-white coloring of the path $P_n$ with exactly $r$ red vertices. Let $Q$ be the longest subpath of $P_n$ whose two leaves are red. Then the restriction of the coloring $c$ to $Q$ is not a symmetric coloring if and only if the original coloring $c$ is an ID-coloring of $P_n$.
\ethm


In this note, we study ID-colorings of cycles, and we establish a criterion to determine whether a red-white coloring of a cycle is an ID-coloring or not when the order (or size) of the cycle is a prime number.
The approach to this criterion is similar to the criterion for paths: we pay attention to symmetries of given red-white colorings of cycles.
Detailed arguments start from Section 2.

A motivation for finding such criteria for various classes of graphs can be described as follows.
Let $G$ be a connected graph and let $H$ be a connected subgraph of $G$, where $H$ has at least three vertices. Let $c$ be a red-white coloring of $G$ such that all red vertices belong to the subgraph $H$. 
Suppose that the restriction of the coloring $c|_H$ is an ID-coloring of $H$ and $d_G(x, y) =d_H(x, y)$ for every two vertices $x$ and $y$ of $H$.
Then the codes of the red vertices are distinct not only in $H$ but also in $G$, as the existence of white vertices outside $H$ does not change the codes of the vertices of $H$. 
It remains then to determine whether the white vertices of $G$ have distinct codes and whether the coloring $c$ is an ID-coloring of $G$. 
This is a more efficient method of determining whether
the red vertices of $H$ have distinct codes or not in $G$. 
This idea was useful in the proof of a theorem in \cite{IDgraph1}, which is the following:
\bthm        \label{cyclered} 
\  For each integer~$n \ge 6$, there is an ID-coloring of~$C_n$ with exactly~$r$ red vertices if and only if  $3 \le r \le n-3$.
Consequently, $\ID(C_n)=3$ for $n \ge 6$. 
\ethm


This is a theorem concerned about cycles. 
However, it was utilized in the proof that the cycle $C_n$ contains a path $P$ such that (1) $P$ contains all red vertices of $C_n$, (2) the restriction of the coloring (used in the proof) to $P$ is an ID-coloring, and (3) $d_{C_n}(x, y) =d_P(x, y)$ for every two vertices $x$ and $y$ of $P$. 
Since the red vertices of $P$ have distinct codes not only in $P$ but also in $C_n$, it remained to determine whether the white vertices of $C_n$ have distinct codes, which made the proof more efficient. 


On the other hand, suppose again that $G$ is a connected graph and $H$ is a connected subgraph of $G$, where $H$ contains at least two  vertices. Let $c$ be a red-white coloring of $G$ such that all red vertices belong to the subgraph $H$. 
Suppose now that the restriction of the coloring $c|_H$ is {\it not} an ID-coloring of $H$ and $d_G(x, y) =d_H(x, y)$ for every two vertices $x$ and $y$ of $H$.
Then the coloring $c$ is {\it not} an ID-coloring of $G$, because there are at least two vertices sharing the same code in $H$, and these vertices share the same code in $G$ as well, as the existence of white vertices outside $H$ does not change the codes of the vertices of $H$.


\vs{2mm}
The following observation obtained in \cite{IDgraph1} will be useful throughout this paper.

\bpr \label{redwhite} 
\ Let~$c$ be a red-white coloring of a connected graph~$G$ where there is at least one vertex of each color. 
If $x$ is a red vertex and $y$ is a white vertex, then $\vd(x) \ne \vd(y)$.
Equivalently, if $\vd(x) = \vd(y)$, then $x$ and $y$ are both red or both white.
\epr 


\section{The Main Theorem}
\label{cycles}

In this main section, we present and prove the main theorem of this paper.
That is, we establish a criterion to determine whether a red-white coloring of a cycle is an ID-coloring or not, when the order (size) of the cycle is a prime number.

\subsection{Symmetric Colorings of Odd Cycles}

First, we introduce a concept that we need to state the main theorem.
Let $n$ be an odd integer at least 3 and let $u$ be a vertex of the cycle $C_n$.
A red-white coloring $c$ of $C_n$ is called 
{\bf symmetric with respect to the vertex $u$} 
if (1) $c$ assigns either red or white to $u$, and (2) $c$ assigns the same color to the two vertices that have distance $d$ from $u$, for each $d$ ($1\le d \le \flo{\frac n2}$).
The vertex $u$ is called the {\bf central vertex} of the coloring $c$,
and two vertices that are equidistant from $u$ (that have the same color) are called {\bf partners} with respect to $u$. 
For convenience, we consider that the partner of the central vertex is the central vertex itself.
Figure \ref{C13-pw} shows an example of a symmetric coloring of $C_{13}$ with respect to the vertex $u$ (which is labeled). Each pair of partners is indicated by a two-way arrow.

\bfi
\scalebox{0.8}{
{\unitlength 0.1in%
\begin{picture}(20.6200,22.6000)(5.5300,-25.6000)%
%
\special{pn 8}%
\special{ar 1580 1560 1000 1000 0.0000000 6.2831853}%
%
\special{sh 1.000}%
\special{ia 1560 560 54 54 0.0000000 6.2831853}%
\special{pn 13}%
\special{ar 1560 560 54 54 0.0000000 6.2831853}%
%
\special{sh 0}%
\special{ia 1034 730 51 51 0.0000000 6.2831853}%
\special{pn 13}%
\special{ar 1034 730 51 51 0.0000000 6.2831853}%
%
\special{sh 0}%
\special{ia 2144 730 51 51 0.0000000 6.2831853}%
\special{pn 13}%
\special{ar 2144 730 51 51 0.0000000 6.2831853}%
%
\special{sh 1.000}%
\special{ia 780 970 54 54 0.0000000 6.2831853}%
\special{pn 13}%
\special{ar 780 970 54 54 0.0000000 6.2831853}%
%
\special{sh 1.000}%
\special{ia 2390 970 54 54 0.0000000 6.2831853}%
\special{pn 13}%
\special{ar 2390 970 54 54 0.0000000 6.2831853}%
%
\special{sh 0}%
\special{ia 614 1330 51 51 0.0000000 6.2831853}%
\special{pn 13}%
\special{ar 614 1330 51 51 0.0000000 6.2831853}%
%
\special{sh 0}%
\special{ia 2564 1330 51 51 0.0000000 6.2831853}%
\special{pn 13}%
\special{ar 2564 1330 51 51 0.0000000 6.2831853}%
%
\special{sh 0}%
\special{ia 2564 1730 51 51 0.0000000 6.2831853}%
\special{pn 13}%
\special{ar 2564 1730 51 51 0.0000000 6.2831853}%
%
\special{sh 0}%
\special{ia 604 1730 51 51 0.0000000 6.2831853}%
\special{pn 13}%
\special{ar 604 1730 51 51 0.0000000 6.2831853}%
%
\special{sh 1.000}%
\special{ia 740 2080 54 54 0.0000000 6.2831853}%
\special{pn 13}%
\special{ar 740 2080 54 54 0.0000000 6.2831853}%
%
\special{sh 1.000}%
\special{ia 2420 2080 54 54 0.0000000 6.2831853}%
\special{pn 13}%
\special{ar 2420 2080 54 54 0.0000000 6.2831853}%
%
\special{sh 0}%
\special{ia 1004 2380 51 51 0.0000000 6.2831853}%
\special{pn 13}%
\special{ar 1004 2380 51 51 0.0000000 6.2831853}%
%
\special{sh 0}%
\special{ia 2154 2380 51 51 0.0000000 6.2831853}%
\special{pn 13}%
\special{ar 2154 2380 51 51 0.0000000 6.2831853}%
\put(14.9000,-4.3000){\makebox(0,0)[lb]{$u$}}%
%
\special{pn 8}%
\special{pa 1680 987}%
\special{pa 2280 987}%
\special{fp}%
\special{sh 1}%
\special{pa 2280 987}%
\special{pa 2213 967}%
\special{pa 2227 987}%
\special{pa 2213 1007}%
\special{pa 2280 987}%
\special{fp}%
%
\special{pn 8}%
\special{pa 1680 987}%
\special{pa 880 987}%
\special{fp}%
\special{sh 1}%
\special{pa 880 987}%
\special{pa 947 1007}%
\special{pa 933 987}%
\special{pa 947 967}%
\special{pa 880 987}%
\special{fp}%
\special{pa 880 987}%
\special{pa 880 987}%
\special{fp}%
%
\special{pn 8}%
\special{pa 1486 2385}%
\special{pa 2041 2385}%
\special{fp}%
\special{sh 1}%
\special{pa 2041 2385}%
\special{pa 1974 2365}%
\special{pa 1988 2385}%
\special{pa 1974 2405}%
\special{pa 2041 2385}%
\special{fp}%
%
\special{pn 8}%
\special{pa 1486 2385}%
\special{pa 1110 2385}%
\special{fp}%
\special{sh 1}%
\special{pa 1110 2385}%
\special{pa 1177 2405}%
\special{pa 1163 2385}%
\special{pa 1177 2365}%
\special{pa 1110 2385}%
\special{fp}%
\special{pa 1110 2385}%
\special{pa 1110 2385}%
\special{fp}%
%
\special{pn 8}%
\special{pa 1710 2077}%
\special{pa 2310 2077}%
\special{fp}%
\special{sh 1}%
\special{pa 2310 2077}%
\special{pa 2243 2057}%
\special{pa 2257 2077}%
\special{pa 2243 2097}%
\special{pa 2310 2077}%
\special{fp}%
%
\special{pn 8}%
\special{pa 1710 2077}%
\special{pa 840 2077}%
\special{fp}%
\special{sh 1}%
\special{pa 840 2077}%
\special{pa 907 2097}%
\special{pa 893 2077}%
\special{pa 907 2057}%
\special{pa 840 2077}%
\special{fp}%
\special{pa 840 2077}%
\special{pa 840 2077}%
\special{fp}%
%
\special{pn 8}%
\special{pa 1496 747}%
\special{pa 2051 747}%
\special{fp}%
\special{sh 1}%
\special{pa 2051 747}%
\special{pa 1984 727}%
\special{pa 1998 747}%
\special{pa 1984 767}%
\special{pa 2051 747}%
\special{fp}%
%
\special{pn 8}%
\special{pa 1496 747}%
\special{pa 1130 747}%
\special{fp}%
\special{sh 1}%
\special{pa 1130 747}%
\special{pa 1197 767}%
\special{pa 1183 747}%
\special{pa 1197 727}%
\special{pa 1130 747}%
\special{fp}%
\special{pa 1130 747}%
\special{pa 1130 747}%
\special{fp}%
%
\special{pn 8}%
\special{pa 1850 1747}%
\special{pa 2450 1747}%
\special{fp}%
\special{sh 1}%
\special{pa 2450 1747}%
\special{pa 2383 1727}%
\special{pa 2397 1747}%
\special{pa 2383 1767}%
\special{pa 2450 1747}%
\special{fp}%
%
\special{pn 8}%
\special{pa 1850 1747}%
\special{pa 690 1747}%
\special{fp}%
\special{sh 1}%
\special{pa 690 1747}%
\special{pa 757 1767}%
\special{pa 743 1747}%
\special{pa 757 1727}%
\special{pa 690 1747}%
\special{fp}%
\special{pa 690 1747}%
\special{pa 690 1747}%
\special{fp}%
%
\special{pn 8}%
\special{pa 1860 1330}%
\special{pa 2460 1330}%
\special{fp}%
\special{sh 1}%
\special{pa 2460 1330}%
\special{pa 2393 1310}%
\special{pa 2407 1330}%
\special{pa 2393 1350}%
\special{pa 2460 1330}%
\special{fp}%
%
\special{pn 8}%
\special{pa 1860 1330}%
\special{pa 720 1330}%
\special{fp}%
\special{sh 1}%
\special{pa 720 1330}%
\special{pa 787 1350}%
\special{pa 773 1330}%
\special{pa 787 1310}%
\special{pa 720 1330}%
\special{fp}%
\special{pa 720 1330}%
\special{pa 720 1330}%
\special{fp}%
\end{picture}}%
}
\capt{A symmetric coloring of $C_{13}$ with respect to the vertex $u$}
\label{C13-pw}
\efi

The following is a basic fact about red-white colorings of $C_3$ and $C_5$.

\bpr \label{prop-C3-C5-pw}
All red-white colorings of $C_3$ and $C_5$ are symmetric colorings with respect to some vertex of them. 
\epr \vs{-3mm}
\pf\ All possible red-white colorings of $C_3$ and $C_5$ are shown in Figure \ref{C3-C5-pw} and they are symmetric colorings with respect to the central vertex $u$, which is labeled for each figure.
\hfill \qed

\bfi
\scalebox{0.8}{\input{C3-C5-pw.tex}}
\capt{All possible red-white colorings of $C_3$ and $C_{5}$}
\label{C3-C5-pw}
\efi

\subsection{The Statement of the Main Theorem}

The following is the main result that we discuss in this paper. 
\bthm \label{prime-id-pw}
Let $n\ge3$ be a prime number.
A red-white coloring of the cycle $C_n$ is an ID-coloring if and only if it is not a symmetric coloring with respect to any vertex of $C_n$.
\ethm

Naturally, Theorem \ref{prime-id-pw} can be split into the following two statements (using contrapositives).
\vs{-3mm}
\bthm \label{prime-id-pw1}
Let $n\ge3$ be a prime number.
A red-white coloring of the cycle $C_n$ is a symmetric coloring with respect to some vertex of $C_n$ if it is not an ID-coloring of $C_n$.
\ethm
\vs{-3mm}
\bpr \label{prime-id-pw2}
Let $n\ge3$ be a prime number.
A red-white coloring of the cycle $C_n$ is not an ID-coloring if it is a symmetric coloring with respect to some vertex of $C_n$.
\epr


\vs{-1mm}
First, we prove Theorem \ref{prime-id-pw1}.
Due to Proposition \ref{prop-C3-C5-pw}, Theorem \ref{prime-id-pw1} is true for $n=3, 5$, so we consider a prime number $n\ge 7$ and $C_n = (u_1, u_2, \cdots, u_n, u_1)$.
For a red-white coloring $c$ of $C_n$,
we suppose that $c$ is not an ID-coloring of $C_n$.
This means that there are at least two vertices of $C_n$ that have the same code. We may assume that such vertices are $u_1$ and $u_k$, where $2\le k \le \cent$. 
Namely, $\vd(u_1)=\vd(u_k)$.
Based on this fact, we will construct a symmetric coloring of $C_n$ with respect to a vertex of $C_n$, using an algorithm that will be introduced later.

Note that we will work on the indices of the vertices of $C_n= (u_1, u_2, \cdots, u_n, u_1)$, and all the computations will be performed in the set $\ZZ/n\ZZ$, where $n=0$.
Since $n$ is a prime number, the set $\ZZ/n\ZZ$ is a field, so all four computations (additions, subtractions, multiplications and divisions) are possible.

\subsection{Tools and Terminologies}

Before we go into the details of a proof of Theorem \ref{prime-id-pw1} (constructing a symmetric coloring of $C_n$), we need some preliminary tools, terminologies and results.
\vs{-1mm}
\bpr \label{central-vertex}
Let $n$ be an odd number with $n\ge 3$ and suppose $C_n=(u_1, \cdots, u_n, u_1)$.
For any integer $k$ with $2 \le k \le \cent$, there is a unique vertex $u_j$ such that $\dist 1j = \dist kj$.
Furthermore,
\vs{-1.5mm}
\[
j = 
\begin{cases}
 \frac{k+1}2 & (k\ is \ odd) \\[2mm]
 \frac{n+k+1}2 & (k\ is\ even).
\end{cases}
\]
\epr \vs{-5mm}
\pf\ \ 
First, suppose that $k$ is odd. 
Then let $j = \frac{k+1}2$ and note that $1<j<k$.
We obtain $\dist 1j = \dist kj = \frac{k-1}2$.
For uniqueness, observe that there is clearly no vertex $\uel$ with $2\le \ell \le k-1$ that is equidistant from $u_1$ and $u_k$ other than $u_j$.
If there is $\ell$ with $k+1 \le \ell \le n$ such that $\dist 1 \ell = \dist k \ell$, 
then the path $(u_k, \cdots, \uel, \cdots, u_n, u_1)$ has an even length $2\cdot\dist 1\ell = 2\cdot\dist k\ell$. 
However, it is impossible given that $C_n$ is an odd cycle and the path $(u_1, \cdots, u_k)$ has an even length, $k-1$.

Next, suppose that $k$ is even.
Then let $j= \frac{n+k+1}2$ and note that $k < j \le n$.
We obtain $\dist 1j = \dist kj = \frac{n-k+1}2$.
For uniqueness, observe that there is clearly no vertex $\uel$ with $k+1 \le \ell \le n$ that is equidistant from $u_1$ and $u_k$ other than $u_j$.
If there is $\ell$ with $2 \le \ell \le k-1$ such that $\dist 1 \ell = \dist k \ell$, 
then the path $(u_1, \cdots, \uel, \cdots, u_k)$ 
must have an even length $2\cdot\dist 1\ell = 2\cdot\dist k\ell$, which is impossible because this path has an odd length $k-1$.
\hfill
\qed

\vs{3mm}

The unique vertex $u_j$ with respect to the two vertices $u_1$ and $u_k$ determined by Proposition \ref{central-vertex} will be the central vertex of the symmetric coloring that we are planning to construct (and the vertices $u_1$ and $u_k$ will be partners of the symmetric coloring with respect to $u_j$).

Once we determine the central vertex $u_j$ of $C_n$ with respect to $u_1$ and $u_k$ ($2 \le k \le \cent$), there are a few more vertices that will be useful with names, if $n\ge 7$. 
The vertices $u_{j-1}$ and $u_{j+1}$ are called the {\bf semi-central vertices} of $C_n$ with respect to $u_1$ and $u_k$. 
The vertices $u_{j_1}$ and $u_{j_2}$, where $j_1 = j + \flo{\frac n2}$ and $j_2 = j + \cei{\frac n2}$, are called the {\bf anti-central vertices} of $C_n$ with respect to $u_1$ and $u_k$. 
If $k=2$, then the anti-central vertices with respect to $u_1$ and $u_2$ are $u_1$ and $u_2$ themselves.
On the other hand, if $k=3$, then the central vertex with respect to $u_1$ and $u_3$ is $u_2$ and the semi-central vertices with respect to $u_1$ and $u_3$ are $u_1$ and $u_3$ themselves.
Note that the anti-central vertices $u_{j_1}$ and $u_{j_2}$ are consecutive in $C_n$, while the semi-central vertices $u_{j-1}$ and $u_{j+1}$ have distance 2.
Also, keep in mind that $\diam (C_n)=\flo{\frac n2}$.

\vs{1mm}
It is convenient to partition $V(C_n)$ into three sets as follows:
\vs{-3mm}
\ben
\item If $k$ is even, then recall that $j= \frac{n+k+1}2$.
Let $I :=\{ \ujt, u_{j_2 + 1}, \cdots, u_{j-1} \}$ and
let $I' := \{ u_{j+1}, u_{j+2}, \cdots, u_n, u_1, u_2, \cdots, \ujo \}$.
Hence we can express $V(C_n) = I \cup I' \cup \{u_j\}$.
\vs{-1mm}
\item If $k$ is odd, then recall that $j= \frac{k+1}2$.
Let $I :=\{ u_{j+1}, u_{j+2}, \cdots, \ujo \}$ and
let $I' := \{ \ujt, u_{j_2 + 1}, \cdots, u_n, u_1, u_2, \cdots, u_{j-1}  \}$.
Hence we can express $V(C_n) = I \cup I' \cup \{u_j\}$.
\een
\vs{-2mm}
In both cases, note that $I$ and $I'$ are sets of consecutive vertices in $C_n$, which start and end with one of the anti-central vertices and one of the semi-central vertices.
Also, it is important to keep in mind that $u_k\in I$, while $u_1, u_n\in I'$, no matter whether $k$ is even or odd.
Note that neither $I$ nor $I'$ contains $u_j$.
For convenience, we sometimes let $I$ and $I'$ denote the sets of indices of the corresponding vertices.
For example, if $k$ is even, then we sometimes let 
$I :=\{ j_2,  j_2 + 1, \cdots, j-1 \}$ and
$I' := \{ j+1, j+2, \cdots, n, 1, 2, \cdots, j_1 \}$.

In order to illustrate the new terminologies and notations that we have introduced, let us consider the cycle $C_{11} = (u_1, u_2, \cdots, u_{11}, u_1)$.
\vs{-3mm}
\ben
\item For the two vertices $u_1$ and $u_4$ ($k=4$), we determine the central vertex with respect to them, which is $u_8$ ($j=8$). 
The vertices $u_7$ and $u_9$ are the semi-central vertices ($j-1=7$ and $j+1=9$), and the vertices $u_2$ and $u_3$ are the anti-central vertices ($j_1 = j+ \flnt = 2$ and $j_2 = j+\cent = 3$).
The sets $I$ and $I'$ represent $\{ u_3, u_4, u_5, u_6, u_7 \}$ and $\{ u_9, u_{10}, u_{11}, u_1, u_2 \}$, respectively, and $V(C_{11}) = I \cup I' \cup \{u_8\}$.
\vs{-2mm}
\item For the two vertices $u_1$ and $u_5$ ($k=5$), we determine the central vertex with respect to them, which is $u_3$ ($j=3$). 
The vertices $u_2$ and $u_4$ are the semi-central vertices ($j-1=2$ and $j+1=4$), and the vertices $u_8$ and $u_9$ are the anti-central vertices ($j_1 = j+ \flnt = 8$ and $j_2 = j+\cent = 9$).
The sets $I$ and $I'$ represent $\{ u_4, u_5, u_6, u_7, u_8 \}$ and $\{ u_9, u_{10}, u_{11}, u_1, u_2 \}$, respectively, and $V(C_{11}) = I \cup I' \cup \{u_3\}$.
\een


For a vertex $u_\ell$ of $C_n$, let $u_{\ell'}$ denote the partner of $u_\ell$.
Namely, the index $\ell'$ means the index of the partner of the vertex $u_\ell$.
For convenience, we sometimes write ``$\ell$ and $\ell'$ are partners'',  meaning that $u_\ell$ and $u_{\ell'}$ are partners.
Unless indicated, we will assume that partners are always with respect to the central vertex $u_j$, which is with respect to the two vertices $u_1$ and $u_k$.
It is important to note the following:
\vs{-2mm}
\bobs \label{IandI'}
If $u_\ell$ and $u_{\ell'}$ are partners that are not $u_j$, then
\vs{-1.5mm}
\[
u_\ell \in I   \iff  u_{\ell'} \in I'.
\]
\eobs
\vs{-3mm}
\bpr \label{partners-index}
The vertices $\uel$ and $\uelp$ are partners if and only if
\vs{-1.5mm}
\[
\ell' = n+k+1-\ell.
\]
\epr \vs{-3mm}
\pf \ \ 
First, suppose that $\ell$ and $\ell'$ are partners.
If $\ell = j$, then the partner of $j$ is $j$ itself, and
the equality actually holds:
if $k$ is even, then $j= \frac{n+k+1}2$,
so
$j'= n+k+1-j=n+k+1- \frac{n+k+1}2=\frac{n+k+1}2=j$.
On the other hand, if $k$ is odd, then $j= \frac{k+1}2$,
so
$j'= n+k+1-j=n+k+1- \frac{k+1}2=n +\frac{k+1}2=\frac{k+1}2=j$.
So we now assume that $\ell \ne j$.
Suppose that $k$ is even (namely $j= \frac{n+k+1}2$ and $k<j \le n$).
If $\ell \in I$, then  $\dist \ell j = j-\ell$,
so $\ell' = j+(j-\ell)= 2j-\ell = n+k+1-\ell$. 
If $\ell \in I'$, then $\dist \ell j = \ell - j$,
so $\ell' = j - (\ell -j) = 2j - \ell = n+k+1 -\ell$.
On the other hand, suppose that $k$ is odd (namely $j= \frac{k+1}2$ and $1<j<k$).
If $\ell \in I$, then $\dist \ell j = \ell - j$,
so $\ell' = j - (\ell - j)= 2j-\ell = n+k+1-\ell$. 
If $\ell \in I'$, then $\dist \ell j = j - \ell$,
so $\ell' = j + (j- \ell ) = 2j - \ell = n+k+1 -\ell$.

For the converse, suppose that $\ell$ and $\ell'$ satisfy $\ell' = n+k+1- \ell$.
First, if $\ell = j$, then we immediately get $\ell = j = \ell'$ no matter whether $k$ is even or odd.
Since the partner of $j$ is $j$ itself, this is a correct result.
So we now assume that $\ell \ne j$.
Suppose that $k$ is even (namely $j= \frac{n+k+1}2$ and $k<j \le n$).
If $\ell \in I$, then $\dist \ell j = j-\ell = \frac{n+k+1}2 - \ell = n+k+1-\ell - \frac{n+k+1}2 = \ell' - j$.
Note that $j+ \dist \ell j = j+ (\ell'-j) = \ell' \in I'$.
This means that $\ell$ and $\ell'$ are distinct and equidistant from $j$,
and hence $\ell$ and $\ell'$ are partners.
On the other hand,  if $\ell \in I'$, then $\dist \ell j =\ell-j = \ell - \frac{n+k+1}2 = \frac{n+k+1}2 - (n+k+1-\ell) = j-\ell'$.
Note that $j - \dist \ell j = j -  (j - \ell') = \ell'\in I$.
This means that $\ell$ and $\ell'$ are distinct and equidistant from $j$,
and hence $\ell$ and $\ell'$ are partners.

Next, suppose that $k$ is odd (namely $j= \frac{k+1}2$ and $1<j<k$).
If $\ell \in I$, then $\dist \ell j = \ell - j =\ell -  \frac{k+1}2 =\frac{k+1}2 - (n+ k+1-\ell)  = j-\ell'$.
Note that $j -  \dist \ell j = j- (j-\ell') = \ell'\in I'$.
This means that $\ell$ and $\ell'$ are distinct and equidistant from $j$,
and hence $\ell$ and $\ell'$ are partners.
On the other hand, if $\ell \in I'$, then $\dist \ell j =j-\ell = \frac{k+1}2-\ell = (n+k+1-\ell) - \frac{k+1}2 = \ell' - j$.
Note that $j + \dist \ell j = j + (\ell' - j) = \ell'\in I$.
This means that $\ell$ and $\ell'$ are distinct and equidistant from $j$,
and hence $\ell$ and $\ell'$ are partners.
\hfill \qed

\vs{3mm}
We have already mentioned that $u_1$ and $u_k$ are partners. 
We can see this fact using Proposition \ref{partners-index} as well.
Indeed, $k'= n+k+1-k=n+1=1$.
Also, the partner of $u_{a'}$ (the partner of $u_a$) is $u_a$.
Indeed, observe that $(a')' = n+k+1-a' = n+k+1-(n+k+1-a)=a$. 

The following proposition is useful before we discuss Proposition \ref{preserve}.
\vs{-1mm}
\bpr \label{partner-lem}
For any integers $\ell$ and $a$, we have
\vs{-2mm}
\[
(\ell + a)' = \ell' - a = a' - \ell.
\]
\epr \vs{-3mm}
\pf \ 
By Proposition \ref{partners-index}, $(\ell + a)' = n+k+1 - (\ell + a) = n+k+1-\ell - a = \ell' - a = a' - \ell$.
\hfill \qed

\vs{2mm}
The distance between two vertices $u_\al$ and $u_\be$ is preserved after taking the partner of both of them.
\vs{-2mm}
\bpr \label{preserve}
For distinct integers $\al$ and $\be$, the following equality holds: $\dist{\al}{\be} = \dist{\al'}{\be'}$.
\epr \vs{-3mm}
\pf \ 
We may assume that $1 \le \al < \be \le n$.
Suppose that $\dist{\al}{\be} =d$.
Note that $d = \be- \al$ or $\al - \be$ (mod $n$).
\\
(1) If $d = \be - \al$, then the $u_\al - u_\be$ geodesic is 
$P=(u_\al, u_{\al+1}, u_{\al+2}, \cdots, u_{\al+d} = u_\be)$, which has length $d$.
Now let $P'$ be the path whose vertices are the partners of the vertices of $P$. 
By Proposition \ref{partner-lem}, observe that
$P' = (u_{\al'}, u_{{(\al+1)}'}, u_{{(\al+2)}'}, \cdots, u_{{(\al+d)}'} = u_{\be'}) = (u_{\al'}, u_{\al'-1}, u_{\al'-2}, \cdots, u_{\al'-d} = u_{\be'})$ and hence $P'$ has length $d$ as well, and $P'$ is the $u_{\al'}-u_{\be'}$ geodesic in $C_n$.
Therefore, $\dist{\al'}{\be'} = d = \dist{\al}{\be}$.
\\
(2) If $d = \al - \be$, then the $u_\be - u_\al$ geodesic is 
$P=(u_\be, u_{\be+1}, u_{\be+2}, \cdots, u_{\be+d} = u_\al) $ (the indices are all mod $n$), which has length $d$.
Now let $P'$ be the path whose vertices are the partners of the vertices of $P$. 
By Proposition \ref{partner-lem}, observe that
$P' = (u_{\be'}, u_{{(\be+1)}'}, u_{{(\be+2)}'}, \cdots, u_{{(\be+d)}'} = u_{\al'}) = (u_{\be'}, u_{\be'-1}, u_{\be'-2}, \cdots, u_{\be'-d} = u_{\al'})$ and hence $P'$ has length $d$ as well, and $P'$ is the $u_{\be'}-u_{\al'}$ geodesic in $C_n$.
Therefore, $\dist{\be'}{\al'} = d = \dist{\be}{\al}$.
\hfill \qed

\vs{3mm}
Given a red-white coloring $c$ of $C_n$ with vertices $u_a$ and $u_b$, we write $u_a \llra u_b$ if $c$ assigns the same color to $u_a$ and $u_b$ (namely, $u_a$ and $u_b$ are both red or both white).

\vs{2mm}
\subsection{An Algorithm and an Example}
Now we turn our attention back to proving Theorem \ref{prime-id-pw1}.
Let $n\ge 7$ be a prime number and let $C_n = (u_1, u_2, \cdots, u_n, u_1)$.
For a red-white coloring $c$ of $C_n$, we suppose that $c$ is not an ID-coloring of $C_n$.
This means that there are at least two distinct vertices of $C_n$ that have the same code. We may assume that such vertices are $u_1$ and $u_k$, where $2\le k \le \cent$. Namely, $\vd(u_1)=\vd(u_k)$.
Given the two vertices $u_1$ and $u_k$, we can determine the vertex $u_j$ such that $u_1$ and $u_k$ are equidistant from $u_j$ by Proposition \ref{central-vertex}.
Therefore, it suffices to show that the coloring $c$ assigns the same color to the two vertices that have distance $d$ from $u_j$, for each $d$ ($1\le d \le \flnt$), which makes $c$ a symmetric coloring with respect to $u_j$.
The following algorithm allows us to do so.
For convenience, we let $k_1 = k$ and hence $k_1'= k' = 1$.


\balg \label{thealgorithm}
\quad
{\bf Step 0.}
Since $\vd(u_1)=\vd(u_k)$, we immediately obtain $u_1 \llra u_k$.

\noindent
{\bf Step 1.}
Let $d_1 = k-1$.
Since $u_1 \llra u_k$ and the $d_1$-th coordinate of $\vd(u_1)$ and $\vd(u_k)$ are the same, we obtain $u_{k+d_1} = u_{2k-1}  \llra u_{1-d_1} =u_{2-k}=u_{n+2-k} $.
Note that $2k-1$ and $n+2-k$ are partners since $(2k-1) + (n+2-k) = k+1$.

\noindent
\ul{Fact 1:} Either $(2k-1)$ or $(n+2-k)$ is in I.

\noindent
Let $k_2$ be the one which is in I and let $k_2'$ be the other one (note that $k_2$ and $k_2'$ are partners).

\noindent
\ul{Fact 2:} \ $k_2, k_2' \notin \{ j, k, 1 \}$.

\noindent
Now we have obtained $u_{k_2} \llra u_{k_2'}$, which is a pair of partners that are not $u_1$ and $u_k$.

\noindent
{\bf Step s.} ($s\ge 2$)

\noindent
\ul{Fact 3:} Exactly one of $\dist 1{k_s}$ and $\dist k{k_s}$ is $d_{s-1}$.

\noindent
Let $d_s \in \{ \dok{s},\ \dkk{s} \} $ be the one that is not $d_{s-1}$ in Fact 3. Since $k_s \notin \{ 1, k\}$ by the previous step, we have $d_s \ne 0$.

\noindent
\ul{Fact 4:} $d_s\notin \{ d_1,\ \cdots,\ d_{s-1} \}$

\noindent
\ul{Fact 5:} Exactly one of $1-d_s$ and $1+d_s$ belongs to the set $\{ k_s, k_s' \}$.

\noindent
Let $D_s \in \{ 1-d_s,\ 1+d_s \}$ be the one that does \ul{not} belong to the set $\{ k_s, k_s' \}$.

\noindent
If $D_s=j$, then we stop the algorithm.

\noindent
If not,

\noindent
\ul{Fact 6:} either $D_s$ or $D_s'$ is in $I$. 

\noindent
Let $k_{s+1} \in \{ D_s, D_s'\}$ be the one in I and let $k'_{s+1}$ be the other one (note that $k_{s+1}$ and $k'_{s+1}$ are partners).

\noindent
\ul{Fact 7:}  \  $k_{s+1}, k_{s+1}' \notin \{ j, k, 1, k_2, k_2', \cdots, k_s, k_s' \}$.

\noindent
\ul{Fact 8:} $u_{k_{s+1}} \llra u_{k'_{s+1}}$.

\noindent
We run and repeat the algorithm until it terminates (when we obtain $D_s=j$ for some $s$).

\ealg

Let us illustrate how the algorithm works with an example.
Let us consider a red-white coloring $c$ of the cycle $C_7 = (u_1, u_2, \cdots, u_7, u_1)$ (namely $n=7$).
Suppose that $c$ is not an ID-coloring of $C_7$ and suppose that $u_1$ and $u_4$ have the same codes (namely $k=4$ and $\vd(u_1) = \vd(u_4)$).
The central vertex with respect to $u_1$ and $u_4$ is $u_6$ (namely $j=6$).
Note that $I= \{ u_3, u_4, u_5 \}$ and $I' = \{ u_7, u_1, u_2  \}$.
First, we immediately obtain $u_1 \llra u_4$ due to the fact that $\vd(u_1) = \vd(u_4)$ and this is {\bf Step 0} of the algorithm.
Then we move on to {\bf Step 1}.
Let $d_1 = k-1 = 4-1 = 3$.
Since we have $u_1 \llra u_4$ and the 3rd coordinate of the codes $\vd(u_1)$ and $\vd(u_4)$ are the same, we obtain 
$u_{k+d_1} = u_7 \llra u_{1-d_1} = u_5$.
Note that $u_5$ and $u_7$ are partners (with respect to the central vertex $u_6$).
Here, $5 \in I$ and $7\in I'$, so we let $k_2 = 5$ and $k_2' = 7$.
Observe that $k_2 = 5 \ne 6, 4, 1$ and $k_2' = 7 \ne 6, 4, 1$.
Thus we have obtained a pair of partners $u_{5}$ and $u_{7}$ with $u_5 \llra u_7$ and they are neither $u_1$ nor $u_4$ (a pair of partners with the same color that we have already obtained).
Now for {\bf Step 2}, observe that $\dok{2} = \dist 15 = 3$ (caution! It is not 4) and $\dkk{2} = \dist 45 =1$.
Since $\dist 15 = 3 = d_1$, we let $d_2 = \dist 45 = 1$.
Now we consider $1-d_2 = 7$ and $1+d_2 = 2$. Since $2 \ne k_2, k_2'$ (namely $2 \ne 5, 7$), we let $D_2 = 2$.
Since $D_2 = 2 \ne j = 6$, we continue to run the algorithm.
Here, $D_2 = 2 \in I'$ and $D_2' = 3 \in I$,
so we let $k_3 = 3$ and $k_3' = 2$ (note that 2 and 3 are partners).
Observe that $k_3 = 3 \ne 6, 4, 1, 5, 7$ and $k_3' = 2 \ne 6, 4, 1, 5, 7$.
By Fact 7, it turns out that $u_3 \llra u_2$, and they are a pair of partners that are neither the central vertex $u_6$ nor the partner vertices having the same color that we have already obtained.
Lastly, {\bf Step 3}.
Observe that $\dok{3} = \dist 13 = 2$ and $\dkk{3} = \dist 43 =1$.
Since $\dist 43 =1 = d_2$, we let $d_3 = \dist 13 = 2$.
Note that $d_3 \ne d_1, d_2$ (namely $2 \ne 3, 1$).
Now we consider $1-d_3 = 6$ and $1+d_3 = 3$. Since $6 \ne k_3, k_3'$ (namely $6 \ne 3, 2$), we let $D_3 = 6$.
Since $D_3 = 6 = j$, we terminate the algorithm.
So far, we have obtained three pairs of partners, each of which have the same color ($u_1 \llra u_4$, $u_5 \llra u_7$ and $u_2 \llra u_3$), 
so the red-white coloring $c$ of $C_7$ is a symmetric coloring with respect to $u_6$.

For the example we just saw above, we chose $k=4$ and we obtained a symmetric coloring of $C_7$, but we obtain the same results for $k=2, 3$ as well, which means that Theorem \ref{prime-id-pw1} is true for $n=7$.

\subsection{How the Algorithm Works}

Now, we explain why the algorithm works for general $n$ (a prime number at least 7) and $k$ ($2 \le k \le \cent$), proving each ``Fact'' stated in the algorithm.

\vs{2mm}
\noindent
\ul{Fact 1:} Either $(2k-1)$ or $(n+2-k)$ is in I.

\vs{-1mm}
It suffices to show that $2k-1 \ne j$.
Assume to the contrary that $2k-1=j$.
If $k$ is odd, then $2k-1 = \frac{k+1}2$, which is equivalent to $k=1$, which is a contradiction.
If $k$ is even, then $2k-1 = \frac{n+k+1}2 \iff 4k-2 = n+k+1 \iff k=1$, which is again a contradiction (note that $n=0$). 

\vs{1mm}
\noindent
\ul{Fact 2:} \ $k_2, k_2' \notin \{ j, k, 1 \}$.

First, we show that $k_2 \notin \{ j, k, 1\}$.
Since $k_2\in I$, we already have $k_2 \notin \{ j, 1\}$.
Thus, it suffices to show that $k_2 \ne k$.
Assume to the contrary that $k_2 = k$.
If $k_2 = 2k-1$, then $k=2k-1 \iff k=1$, which is a contradiction.
If $k_2 = n+2-k$, then $k=n+2-k \iff k=1$, which is a contradiction as well.
Next, we show that $k_2' \notin \{ j, k, 1 \}$.
Since $k_2'\in I'$, we already have $k_2' \notin \{ j, k\}$.
We also have $k_2' \ne 1$, otherwise $k_2' =1 \iff k_2 = 1' = k$, which is a contradiction. 

\vs{2mm}
\noindent
\ul{Fact 3:} Exactly one of $\dist 1{k_s}$ and $\dist k{k_s}$ is $d_{s-1}$.

If $k_s = 1-d_{s-1}$ or $1+ d_{s-1}$, then $\dok{s}=d_{s-1}$. If $\dkk{s} = d_{s-1}$ as well, then $k_s = j$ by the uniqueness of the central vertex with respect to $u_1$ and $u_k$. However, this is a contradiction to Fact 2 (if $s=2$) and Fact 4 of the Step $s-1$ (if $s\ge 3$). Hence, $\dkk{s}\ne d_{s-1}$.
On the other hand, if $k_s' = 1-d_{s-1}$ or $1+ d_{s-1}$, then $\dokp{s} = \dkk{s} = d_{s-1}$. If $\dok{s}=d_{s-1}$ as well, then $k_s = j$ again, a contradiction. Hence, $\dok{s}\ne d_{s-1}$.

\vs{3mm}
\noindent
\ul{Fact 4:} \ $d_s \notin \{ d_1,\ \cdots,\ d_{s-1} \}$.

By assumption, we already have $d_s \ne d_{s-1}$.
Thus, it suffices to show that $d_s \notin \{ d_1, \cdots, d_{s-2} \}$.
Assume to the contrary that $d_s = d_p$ for some $1 \le p \le s-2$.
If $d_s = \dok{s} = d_p$, then $k_s = 1-d_p$ or $1+d_p$, and hence $k_s \in \{ k_p,\ k_p',\ k_{p+1},\ k_{p+1}' \}$, which is a contradiction by Fact 2 or Fact 7 of the previous step.
On the other hand, if $d_s = \dkk{s} = \dokp{s} = d_p$, then $k_s' = 1-d_p$ or $1+d_p$, and hence $k_s' \in \{ k_p,\ k_p',\ k_{p+1},\ k_{p+1}' \}$, which is again a contradiction.

\vs{3mm}
\noindent
\ul{Fact 5:} Exactly one of $1-d_s$ and $1+d_s$ belong to the set $\{ k_s, k_s' \}$.

(1) Suppose that $\dok{s} =d_s$ and $\dkk{s}=d_{s-1}$.
Then $k_s = 1-d_s$ or $1+d_s$.
\\
(1-i) Suppose $k_s = 1-d_s$.
Since $C_n$ is an odd cycle, $1+d_s \ne k_s$.
Now, we also have $1+d_s \ne k_s'$. 
Indeed, if $1+d_s = k_s'$, then $d_s = \dokp{s} = \dkk{s} = d_{s-1}$, which is a contradiction.
\\
(1-ii) Suppose $k_s = 1+d_s$.
Since $C_n$ is an odd cycle, $1-d_s \ne k_s$.
Now, we also have $1-d_s \ne k_s'$. 
Indeed, if $1-d_s = k_s'$, then $d_s = \dokp{s} = \dkk{s} = d_{s-1}$, which is a contradiction.

(2) Suppose that $\dok{s} =d_{s-1}$ and $\dkk{s}=d_s=\dokp{s}$.
Then $k_s' = 1-d_s$ or $1+d_s$.
\\
(2-i) Suppose $k_s' = 1-d_s$.
Since $C_n$ is an odd cycle, $1+d_s \ne k_s'$.
Now, we also have $1+d_s \ne k_s$. 
Indeed, if $1+d_s = k_s$, then $d_s = \dok{s} = d_{s-1}$, which is a contradiction.
\\
(2-ii) Suppose $k_s' = 1+d_s$.
Since $C_n$ is an odd cycle, $1-d_s \ne k_s'$.
Now, we also have $1-d_s \ne k_s$. 
Indeed, if $1-d_s = k_s$, then $d_s = \dok{s} = d_{s-1}$, which is a contradiction.

\vs{3mm}
\noindent
\ul{Fact 6:} either $D_s$ or $D_s'$ is in $I$.

This immediately follows from Observation \ref{IandI'}, given that $D_s \ne j$.

\vs{2mm}
\noindent
\ul{Fact 7:} \ $k_3, k_3' \notin \{ j, k, 1, k_2, k_2' \}$. More generally, $k_{s+1}, k_{s+1}' \notin \{ j, k, 1, k_2, k_2', \cdots, k_s, k_s' \}$ for $s\ge 3$.

By assumption, $k_{3} \notin \{k_2, k_2'\}$.
Since $k_{3}\in I$, we also have $k_{3} \notin \{ j, 1 \}$.
Thus, it suffices to show that $k_3 \ne k$.
Assume to the contrary that $k_3 = k$.
If $k_3 = 1-d_2$ or $1+d_2$, then $d_2 = \dok{3} = \dok{} = k-1 = d_1$, which is a contradiction.
If $k_3 = k-d_2$ or $k+d_2$, then $0 \ne d_2 = \dkk{3} = \dkk{} = 0$, a contradiction. Therefore, $k_3 \notin \{ j, k, 1, k_2, k_2' \}$.
Now, $k_3' \notin \{k_2, k_2'\}$ by assumption. 
Since $k_3'\in I'$, we also have $k_3' \notin \{j, k\}$. 
We have $k_3' \ne 1$ as well, otherwise $k_3 = (k_3')'=1'=k$, which is a contradiction by the argument above.
Therefore, $k_3, k_3' \notin \{ j, k, 1, k_2, k_2' \}$.

More generally, let $s \ge 3$.
By assumption, $k_{s+1} \notin \{k_s, k_s'\}$.
Since $k_{s+1}\in I$, we also have $k_{s+1} \notin \{ j, 1, k_2', \cdots, k_{s-1}'\}$.
Thus, we need to show that $k_{s+1} \notin \{ k, k_2, \cdots, k_{s-1}\}$.
Note that $k_{s+1} \in \{ 1-d_s, 1+d_s, (1-d_s)', (1+d_s)' \} = \{ 1-d_s, 1+d_s, k-d_s, k+d_s \}$ (observe that $(1-d_s)' = n+k+1- (1-d_s) = k+d_s$ and $(1+d_s)' = n+k+1- (1+d_s) = k-d_s$).
Also, note that $d_s \notin \{d_1, \cdots, d_{s-1} \}$ by Fact 4.
\ben
\item If $k_{s+1} = 1-d_s$ or $1+d_s$, then $\dok{s+1}=d_s$.
If $k_{s+1} = k$, then $d_s = \dok{s+1} = d(u_1, u_k) = k-1 = d_1$, which is a contradiction.
On the other hand, if $k_{s+1} = k_p$ for $2 \le p \le s-1$,
then $d_s = \dok{s+1} = \dok{p} \in \{ d_{p-1}, d_p\}$, which is also a contradiction.
\item If $k_{s+1} = k-d_s$ or $k+d_s$, then $\dkk{s+1} = d_s$.
If $k_{s+1} = k$, then $0 \ne d_s = \dkk{s+1} = d(u_k, u_k) = 0$, which is a contradiction.
On the other hand, if $k_{s+1} = k_p$ for $2 \le p \le s-1$,
then $d_s = \dkk{s+1} = \dkk{p} \in \{ d_{p-1}, d_p\}$, which is also a contradiction.
\een
Therefore, $k_{s+1} \notin \{ j, k, 1, k_2, k_2', \cdots, k_s, k_s' \}$.
Lastly, note that $k_{s+1}' \notin \{k_s, k_s' \}$ by assumption.
Since $k_{s+1}' \in I'$, we also have $k_{s+1}' \notin \{ j, k, k_2, \cdots, k_{s-1}\}$.
Observe that $k_{s+1}' \notin \{ 1, k_2', \cdots, k_{s-1}' \}$ as well, otherwise $k_{s+1} \in \{ k, k_2, \cdots, k_{s-1} \}$, which is a contradiction by the previous argument.
Therefore, $k_{s+1}' \notin \{ j, k, 1, k_2, k_2', \cdots, k_s, k_s' \}$.

\vs{3mm}
Before we prove Fact 8, it is convenient to establish the following rules.
\vs{-3mm}
\bpr \label{8-nice-rules}
For Step $s$ of the algorithm ($s \ge 2$),
we have the following rules.

(I-1)
If $1-d_s= k_{s+1}$ and $d_s = \dok{s}$, then $1+d_s = k_s$, $k-d_s = k_s'$ and $k+d_s = k_{s+1}'$.

(I-2)
If $1-d_s= k_{s+1}$ and $d_s = \dkk{s}$, then $1+d_s = k_s'$, $k-d_s = k_s$ and $k+d_s = k_{s+1}'$.

(II-1)
If $1+d_s= k_{s+1}$ and $d_s = \dok{s}$, then $1-d_s = k_s$,  $k-d_s = k_{s+1}'$ and $k+d_s = k_s'$.

(II-2)
If $1+d_s= k_{s+1}$ and $d_s = \dkk{s}$, then $1-d_s = k_s'$, $k-d_s = k_{s+1}'$ and $k+d_s = k_s$.

(III-1)
If $1-d_s= k_{s+1}'$ and $d_s = \dok{s}$, then $1+d_s = k_s$, $k-d_s = k_s'$ and $k+d_s = k_{s+1}$.

(III-2)
If $1-d_s= k_{s+1}'$ and $d_s = \dkk{s}$, then $1+d_s = k_s'$, $k-d_s = k_s$ and $k+d_s = k_{s+1}$.

(IV-1)
If $1+d_s= k_{s+1}'$ and $d_s = \dok{s}$, then $1-d_s = k_s$, $k-d_s = k_{s+1}$ and $k+d_s = k_s'$.

(IV-2)
If $1+d_s= k_{s+1}'$ and $d_s = \dkk{s}$, then $1-d_s = k_s'$, $k-d_s = k_{s+1}$ and $k+d_s = k_s$.
\epr \vs{-1mm}
\pf \ 
(I-1)
If $1-d_s= k_{s+1}$ and $d_s = \dok{s}$, then
$k_s = 1-d_s$ or $1+d_s$.
Since $1-d_s = k_{s+1}$, it follows that $k_s = 1+d_s$.
Now $k_s' = (1+d_s)' = n+k+1-(1+d_s)=n+k-d_s = k-d_s$.
Also, observe that
$k_{s+1}'=(1-d_s)'=n+k+1-(1-d_s)=n+k+d_s=k+d_s$.

(I-2)
If $1-d_s= k_{s+1}$ and $d_s = \dkk{s} = \dokp{s}$, then
$k_s' = 1-d_s$ or $1+d_s$.
Since $1-d_s = k_{s+1}$, it follows that $k_s' = 1+d_s$.
Now $k_s = (1+d_s)' = n+k+1-(1+d_s)=n+k-d_s = k-d_s$.
Also, observe that
$k_{s+1}'=(1-d_s)'=n+k+1-(1-d_s)=n+k+d_s=k+d_s$.

(II-1)
If $1+d_s= k_{s+1}$ and $d_s = \dok{s}$, then
$k_s = 1-d_s$ or $1+d_s$.
Since $1+d_s = k_{s+1}$, it follows that $k_s= 1-d_s$.
Now $k_s' = (1-d_s)' = n+k+1-(1-d_s)=n+k+d_s = k+d_s$.
Also, observe that
$k_{s+1}'=(1+d_s)'=n+k+1-(1+d_s)=n+k-d_s=k-d_s$.

(II-2)
If $1+d_s= k_{s+1}$ and $d_s = \dkk{s} = \dokp{s}$, then
$k_s' = 1-d_s$ or $1+d_s$.
Since $1+d_s = k_{s+1}$, it follows that $k_s' = 1-d_s$.
Now $k_s = (1-d_s)' = n+k+1-(1-d_s)=n+k+d_s = k+d_s$.
Also, observe that
$k_{s+1}'=(1+d_s)'=n+k+1-(1+d_s)=n+k-d_s=k-d_s$.

(III-1)
If $1-d_s= k_{s+1}'$ and $d_s = \dok{s}$, then
$k_s = 1-d_s$ or $1+d_s$.
Since $1-d_s= k_{s+1}'$, it follows that $k_s = 1+d_s$.
Now $k_s' = (1+d_s)' = n+k+1-(1+d_s)=n+k-d_s = k-d_s$.
Also, observe that
$k_{s+1}=(1-d_s)'=n+k+1-(1-d_s)=n+k+d_s=k+d_s$.

(III-2)
If $1-d_s= k_{s+1}'$ and $d_s = \dkk{s} = \dokp{s}$, then
$k_s' = 1-d_s$ or $1+d_s$.
Since $1-d_s= k_{s+1}'$, it follows that $k_s' = 1+d_s$.
Now $k_s = (1+d_s)' = n+k+1-(1+d_s)=n+k-d_s = k-d_s$.
Also, observe that
$k_{s+1}=(1-d_s)'=n+k+1-(1-d_s)=n+k+d_s=k+d_s$.

(IV-1)
If $1+d_s= k_{s+1}'$ and $d_s = \dok{s}$, then
$k_s = 1-d_s$ or $1+d_s$.
Since $1+d_s= k_{s+1}'$, it follows that $k_s = 1-d_s$.
Now $k_s' = (1-d_s)' = n+k+1-(1-d_s)=n+k+d_s = k+d_s$.
Also, observe that
$k_{s+1}=(1+d_s)'=n+k+1-(1+d_s)=n+k+d_s=k-d_s$.

(IV-2)
If $1+d_s= k_{s+1}'$ and $d_s = \dkk{s} = \dokp{s}$, then
$k_s' = 1-d_s$ or $1+d_s$.
Since $1+d_s= k_{s+1}'$, it follows that $k_s' = 1-d_s$.
Now $k_s = (1-d_s)' = n+k+1-(1-d_s)=n+k+d_s = k+d_s$.
Also, observe that
$k_{s+1}=(1+d_s)'=n+k+1-(1+d_s)=n+k-d_s=k-d_s$.
\hfill\qed

\vs{4mm}
Now, we are prepared to prove Fact 8.

\vs{2mm}
\noindent
\ul{Fact 8:} $u_{k_{s+1}} \llra u_{k_{s+1}'}$.

Let us consider the red-white coloring $c$ of $C_n$ as the following function:

\vs{-8mm}
\[
c(v) = 
\begin{cases}
 1 & (v\ is \ red) \\[2mm]
 0 & (v\ is\ white).
\end{cases}
\]

\vs{-5mm}
Note that $c(u_{k_s}) = c(u_{k_s'})$ since $u_{k_s} \llra u_{k_s'}$ from the previous step.
With the aid of Proposition \ref{8-nice-rules}, we consider the following cases.

(I-1)
If $1-d_s= k_{s+1}$ and $d_s = \dok{s}$, then $1+d_s = k_s$, $k-d_s = k_s'$ and $k+d_s = k_{s+1}'$.
Now observe that the $d_s$-th coordinate of the code $\vd(u_1)$ is
$c(u_{1-d_s}) + c(u_{1+d_s}) = c(u_{k_{s+1}}) + c(u_{k_s})$,
while the $d_s$-th coordinate of the code $\vd(u_k)$ is
$c(u_{k-d_s}) + c(u_{k+d_s}) = c(u_{k_s'}) + c(u_{k_{s+1}'})$. It follows 
\\[-4.5mm]
that $c(u_{k_{s+1}}) = c(u_{k_{s+1}'})$ and hence $u_{k_{s+1}} \llra u_{k_{s+1}'}$, given that $\vd(u_1) = \vd(u_k)$.

(I-2)
If $1-d_s= k_{s+1}$ and $d_s = \dkk{s}$, then $1+d_s = k_s'$, $k-d_s = k_s$ and $k+d_s = k_{{s+1}}'$.
Now observe that the $d_s$-th coordinate of the code $\vd(u_1)$ is
$c(u_{1-d_s}) + c(u_{1+d_s}) = c(u_{k_{s+1}}) + c(u_{k_s'})$,
while the $d_s$-th coordinate of the code $\vd(u_k)$ is
$c(u_{k-d_s}) + c(u_{k+d_s}) = c(u_{k_s}) + c(u_{k_{s+1}'})$.
It follows  \\[-4.5mm]
that $c(u_{k_{s+1}}) = c(u_{k_{s+1}'})$ and hence $u_{k_{s+1}} \llra u_{k_{s+1}'}$, given that $\vd(u_1) = \vd(u_k)$.

(II-1)
If $1+d_s= k_{s+1}$ and $d_s = \dok{s}$, then $1-d_s = k_s$, $k-d_s = k_{{s+1}}'$ and $k+d_s = k_s'$.
Now observe that the $d_s$-th coordinate of the code $\vd(u_1)$ is
$c(u_{1-d_s}) + c(u_{1+d_s}) = c(u_{k_s}) + c(u_{k_{s+1}})$,
while the $d_s$-th coordinate of the code $\vd(u_k)$ is
$c(u_{k-d_s}) + c(u_{k+d_s}) = c(u_{k_{s+1}'}) + c(u_{k_s'})$.
It follows \\[-4.5mm]
that $c(u_{k_{s+1}}) = c(u_{k_{s+1}'})$ and hence $u_{k_{s+1}} \llra u_{k_{s+1}'}$, given that $\vd(u_1) = \vd(u_k)$.

(II-2)
If $1+d_s = k_{s+1}$ and $d_s = \dkk{s}$, then $1-d_s = k_s'$,  $k-d_s = k_{{s+1}}'$ and $k+d_s = k_s$.
Now observe that the $d_s$-th coordinate of the code $\vd(u_1)$ is
$c(u_{1-d_s}) + c(u_{1+d_s}) = c(u_{k_s'}) + c(u_{k_{s+1}})$,
while the $d_s$-th coordinate of the code $\vd(u_k)$ is
$c(u_{k-d_s}) + c(u_{k+d_s}) = c(u_{k_{s+1}'}) + c(u_{k_s})$.
It follows \\[-4.5mm]
that $c(u_{k_{s+1}}) = c(u_{k_{s+1}'})$ and hence $u_{k_{s+1}} \llra u_{k_{s+1}'}$, given that $\vd(u_1) = \vd(u_k)$.

(III-1)
If $1-d_s= k_{s+1}'$ and $d_s = \dok{s}$, then $1+d_s = k_s$, $k-d_s = k_s'$ and $k+d_s = k_{s+1}$.
Now observe that the $d_s$-th coordinate of the code $\vd(u_1)$ is
$c(u_{1-d_s}) + c(u_{1+d_s}) = c(u_{k_{s+1}'}) + c(u_{k_s})$,\\[-4.5mm]
while the $d_s$-th coordinate of the code $\vd(u_k)$ is
$c(u_{k-d_s}) + c(u_{k+d_s}) = c(u_{k_s'}) + c(u_{k_s+1})$.
It follows that $c(u_{k_{s+1}}) = c(u_{k_{s+1}'})$ and hence $u_{k_{s+1}} \llra u_{k_{s+1}'}$, given that $\vd(u_1) = \vd(u_k)$.

(III-2)
If $1-d_s= k_{s+1}'$ and $d_s = \dkk{s}$, then $1+d_s = k_s'$, $k-d_s = k_s$ and $k+d_s = k_{s+1}$.
Now observe that the $d_s$-th coordinate of the code $\vd(u_1)$ is
$c(u_{1-d_s}) + c(u_{1+d_s}) = c(u_{k_s+1'}) + c(u_{k_s'})$,
while the $d_s$-th coordinate of the code $\vd(u_k)$ is
$c(u_{k-d_s}) + c(u_{k+d_s}) = c(u_{k_s}) + c(u_{k_s+1})$.
It follows that $c(u_{k_{s+1}}) = c(u_{k_{s+1}'})$ and hence $u_{k_{s+1}} \llra u_{k_{s+1}'}$, given that $\vd(u_1) = \vd(u_k)$.

(IV-1)
If $1+d_s= k_{s+1}'$ and $d_s = \dok{s}$, then $1-d_s = k_s$, $k-d_s = k_{s+1}$ and $k+d_s = k_s'$.
Now observe that the $d_s$-th coordinate of the code $\vd(u_1)$ is
$c(u_{1-d_s}) + c(u_{1+d_s}) = c(u_{k_s}) + c(u_{k_s+1'})$,
while the $d_s$-th coordinate of the code $\vd(u_k)$ is
$c(u_{k-d_s}) + c(u_{k+d_s}) = c(u_{k_s+1}) + c(u_{k_s'})$.
It follows that $c(u_{k_{s+1}}) = c(u_{k_{s+1}'})$ and hence $u_{k_{s+1}} \llra u_{k_{s+1}'}$, given that $\vd(u_1) = \vd(u_k)$.

(IV-2)
If $1+d_s= k_{s+1}'$ and $d_s = \dkk{s}$, then $1-d_s = k_s'$, $k-d_s = k_{s+1}$ and $k+d_s = k_s$.
Now observe that the $d_s$-th coordinate of the code $\vd(u_1)$ is
$c(u_{1-d_s}) + c(u_{1+d_s}) = c(u_{k_s'}) + c(u_{k_s+1'})$,
while the $d_s$-th coordinate of the code $\vd(u_k)$ is
$c(u_{k-d_s}) + c(u_{k+d_s}) = c(u_{k_s+1}) + c(u_{k_s})$.
It follows that $c(u_{k_{s+1}}) = c(u_{k_{s+1}'})$ and hence $u_{k_{s+1}} \llra u_{k_{s+1}'}$, given that $\vd(u_1) = \vd(u_k)$.


\vs{6mm}

So far, we have shown all the Facts of the algorithm.
However, in order to obtain a symmetric coloring of $C_n$ with respect to the vertex $u_j$, we need to make sure that the algorithm runs by the end of Step $\frac{n-3}2$, 
which is when we verify that the $(\frac{n-1}2)$-th pair of partners obtain the same color.
This means that we need to verify that $D_s \ne j$ for all $2 \le s \le \frac{n-3}2$.
Recall that $D_s$ is either $1-d_s$ or $1+d_s$ (that is neither $k_s$ nor $k_s'$).
If $D_s = 1-d_s$, then $D_s\ne j  \iff  d_s \ne 1-j$.
If $D_s = 1+d_s$, then $D_s\ne j  \iff  d_s \ne j-1$.
Proposition \ref{dsnot1j} guarantees that these are true for all $2 \le s \le \frac{n-3}2$.

\vs{-2mm}
\ble \label{dssks}
$d_s = sk-s$ or $s-sk$ for $s\ge 1$.
\ele \vs{-3mm}
\pf\ 
By definition, we have $d_1 = k-1$.
Recall that $d_2$ is either $\dok{2}$ or $\dkk{2}$ that is not $d_1$.
Also, recall that $k_2 = 2k-1$ or $n+2-k$.\\
(I) Suppose that $k_2=2k-1$. This means that $\dkk{2}=k-1=d_1$, so $\dok{2}=d_2$. Thus $k_2 = 1-d_2$ or $1+d_2$.
If $k_2 = 1-d_2$, then $2k-1=1-d_2 \iff d_2=2-2k$. 
On the other hand, if $k_2 = 1+d_2$, then $2k-1=1+d_2 \iff d_2=2k-2$. \\
(II) Suppose that $k_2=n+2-k=(2k-1)'$. 
Then $k_2' = 2k-1$ and $\dok{2}=\dkkp{2}=k-1=d_1$, so $\dkk{2}=d_2$. 
Thus $k_2 = k-d_2$ or $k+d_2$.
If $k_2 = k-d_2$, then $n+2-k=k-d_2 \iff d_2=2k-2$. 
On the other hand, if $k_2 = k+d_2$, then $n+2-k=k+d_2 \iff d_2=2-2k$.

We have verified the statement for $s=1, 2$.
Now we proceed by induction on $s$.
Suppose that $d_t = tk-t$ or $t-tk$ for all $t = 1, \cdots, s$ for some $s \ge 2$.
We will show that $d_{t+1}=(t+1)k-(t+1)$ or $(t+1)-(t+1)k$.
Notice that $d_{t+1} = \dok{t+1}$ or $\dkk{t+1}$ and $d_{t+1} \notin \{d_1, \cdots, d_t \}$.
We have $k_{t+1} = 1-d_{t+1},\ 1+d_{t+1},\ k-d_{t+1}$  or $k+d_{t+1}$.

(A) Suppose that $k_{t+1} = 1-d_{t+1}$.
Then $d_{t+1} = \dok{t+1}$ and $\dkk{t+1}=d_t$.
Thus $k_{t+1} = k-d_t$ or $k+d_t$.\\
(A-1) Suppose $k_{t+1} = k-d_t$.
Then $1-d_{t+1} = k-d_t \iff d_{t+1}=1-k+d_t$.
If $d_t= t-tk$, then $d_{t+1} = 1-k+(t-tk) = (t+1)-(t+1)k$, which is desired.
On the other hand, if $d_t= tk-t$, then $d_{t+1} = 1-k+(tk-t) = (t-1)k -(t-1)$. 
Since $d_{t+1} \ne d_{t-1}$, it follows that
$d_{t-1} = (t-1) - (t-1)k$. Now $d_{t-1} + d_{t+1} = 0$.
However, observe that $1 \le d_{t-1} \le \flnt$ and $1 \le d_{t+1} \le \flnt$, so $2 \le d_{t-1} + d_{t+1} \le n-1$, which is a contradiction.\\
(A-2) Suppose $k_{t+1} = k+d_t$.
Then $1-d_{t+1} = k+d_t \iff d_{t+1}=1-k-d_t$.
If $d_t= tk-t$, then $d_{t+1} = 1-k-(tk-t) = (t+1)-(t+1)k$, which is desired.
On the other hand, if $d_t = t-tk$, then $d_{t+1} = 1-k-(t-tk) = (t-1)k -(t-1)$. 
This leads to the contradiction same as above.

(B) Suppose that $k_{t+1} = 1+d_{t+1}$.
Then $d_{t+1} = \dok{t+1}$ and $\dkk{t+1}=d_t$.
Thus $k_{t+1} = k-d_t$ or $k+d_t$.\\
(B-1) Suppose $k_{t+1} = k-d_t$.
Then $1+d_{t+1} = k-d_t \iff d_{t+1}=k-1-d_t$.
If $d_t= t-tk$, then $d_{t+1} = k-1-(t-tk) = (t+1)k-(t+1)$, which is desired.
On the other hand, if $d_t= tk-t$, then $d_{t+1} = k-1-(tk-t) = (t-1) -(t-1)k$. 
This leads to the contradiction same as above.\\
(B-2) Suppose $k_{t+1} = k+d_t$.
Then $1+d_{t+1} = k+d_t \iff d_{t+1}=k-1+d_t$.
If $d_t= tk-t$, then $d_{t+1} = k-1+(tk-t) = (t+1)k-(t+1)$, which is desired.
On the other hand, if $d_t= t-tk$, then $d_{t+1} = k-1+(t-tk) = (t-1) -(t-1)k$. 
This leads to the contradiction same as above.

(C) Suppose that $k_{t+1} = k-d_{t+1}$.
Then $d_{t+1} = \dkk{t+1}$ and $\dok{t+1}=d_t$.
Thus $k_{t+1} = 1-d_t$ or $1+d_t$.
If $k_{t+1} = 1-d_t$, then $k-d_{t+1} = 1-d_t \iff d_{t+1}=k-1+d_t$, which is the same as (B-2). 
On the other hand, if $k_{t+1} = 1+d_t$, then $k-d_{t+1} = 1+d_t \iff d_{t+1}=k-1-d_t$, which is the same as (B-1).

(D) Suppose that $k_{t+1} = k+d_{t+1}$.
Then $d_{t+1} = \dkk{t+1}$ and $\dok{t+1}=d_t$.
Thus $k_{t+1} = 1-d_t$ or $1+d_t$.
If $k_{t+1} = 1-d_t$, then $k+d_{t+1} = 1-d_t \iff d_{t+1}=1-k-d_t$, which is the same as (A-2). 
On the other hand, if $k_{t+1} = 1+d_t$, then $k+d_{t+1} = 1+d_t \iff d_{t+1}=1-k+d_t$, which is the same as (A-1).
\hfill \qed

\bpr \label{dsnot1j}
Let $n \ge 7$ be a prime number.
If $1 \le s \le \frac{n-3}2$, then $d_s \notin \{1-j, j-1\}$.
\epr \vs{-3mm}
\pf\ 
Observe that $1 \le s \le \frac{n-3}2  \iff  2 \le 2s \le n-3  \iff 3 \le 2s+1 \le n-2  \iff 1 \le 2s-1 \le n-4$.
In particular, $2s-1$ and $2s+1$ are not 0.
Note that $j = \frac{k+1}2 = \frac{n+k+1}2$ in the set $\ZZ/n\ZZ$ where $n=0$.
Furthermore, since $n$ is a prime number, the set $\ZZ/n\ZZ$ is a field, where division of its elements is well-defined.
By Lemma \ref{dssks}, $d_s = sk-s$ or $s-sk$ for $s\ge 1$.
\\
(A) Suppose $d_s = sk-s$. \\
(A-1) If $d_s = 1-j = 1 - \frac{k+1}2  = \frac{1-k}2$, then $sk-s=\frac{1-k}2 \iff 2sk-2s = 1-k  \iff  (2s+1)k=2s+1 \iff k=1$ (since $2s+1 \ne 0$) and this is a contradiction.\\
(A-2) If $d_s = j-1= \frac{k+1}2 - 1 = \frac{k-1}2$, then $sk-s=\frac{k-1}2   
\iff  2sk -2s = k-1 \iff (2s-1)k=2s-1  \iff k=1$ (since $2s-1 \ne 0$) and this is a contradiction.\\
(B) Suppose $d_s = s-sk$. \\
(B-1) If $d_s=1-j=\frac{1-k}2$, then $s-sk=\frac{1-k}2  \iff  sk-s = \frac{k-1}2$
and this is the same as (A-2).\\
(B-2) If $d_s=j-1=\frac{k-1}2$, then $s-sk=\frac{k-1}2  \iff  sk-s = \frac{1-k}2$
and this is the same as (A-1).
\hfill \qed

\vs{3mm}

With the aid of Proposition \ref{dsnot1j},
we can see that the algorithm runs sufficiently many times, until the end of Step $\frac{n-3}2$, if $n\ge 7$ is a prime number, and we obtain all the pairs of partners having the same color.
Now the proof of Theorem \ref{prime-id-pw1} is completed.

\vs{2mm}
\subsection{The Converse of Theorem \ref{prime-id-pw1} }

Now, it remains to prove Proposition \ref{prime-id-pw2} (the converse of Theorem \ref{prime-id-pw1}) to complete the proof of the main theorem, Theorem \ref{prime-id-pw}.
Let us re-state Proposition \ref{prime-id-pw2} here.
\\[3mm]
{\bf Proposition \ref{prime-id-pw2} } \ 
Let $n\ge3$ be a prime number.
A red-white coloring of the cycle $C_n$ is not an ID-coloring if it is a symmetric coloring with respect to some vertex of $C_n$.

\vs{2mm}
The following proposition will be the essential key to prove Proposition \ref{prime-id-pw2}.

\vs{-2mm}
\bpr \label{pw-code}
A red-white coloring of an odd cycle $C_n$ is a symmetric coloring with respect to some vertex of $C_n$ if and only if each pair of partners of $C_n$ have the same code.
\epr
\vs{-3mm}
\pf \quad  
Suppose that each pair of partners of $C_n$ with respect to $u$ have the same code. Then each pair of partners of $C_n$ must have the same color, since having the same code implies having the same color. Therefore, the coloring $c$ is a symmetric coloring of $C_n$.

Now we suppose that $c$ is a symmetric coloring of $C_n$.
Let $C_n = (u_1, u_2, \cdots, u_n, u_1)$ and $u=u_1$ (i.e. $u_1$ is the central vertex).
Let us consider $c$ as a function on $V(C_n)$ where $c(v)=1$ if $v$ is red and $c(v)=0$  if $v$ is white.
For $1 \le d \le \flnt$, observe that 
\begin{eqnarray*}
\vd(u_{1+d}) = & (\ c(u_{1+d-1}) + c(u_{1+d+1}), \ c(u_{1+d-2}) + c(u_{1+d+2}),\ \cdots, \\
 & c(u_{1+d-a}) + c(u_{1+d+a}), \cdots, c(u_{1+d-\flnt}) + c(u_{1+d+\flnt}) \ ), \\[-4mm]
 \vd(u_{1-d}) = & (\ c(u_{1-d-1}) + c(u_{1-d+1}), \ c(u_{1-d-2}) + c(u_{1-d+2}),\ \cdots, \\
 & c(u_{1-d-a}) + c(u_{1-d+a}), \cdots, c(u_{1-d-\flnt}) + c(u_{1-d+\flnt}) \ ).
\end{eqnarray*}
\vs{-6mm}

Now, for $1 \le a \le \flnt$, notice that $u_{1+d+a}$ and $u_{1-d-a}$ are partners (since $1-d-a=1-(d+a)$), and $u_{1+d-a}$ and $u_{1-d+a}$ are partners (since $1+d-a = 1+(d-a)$ and $1-d+a = 1-(d-a)$).
Therefore, for any $1 \le d \le \flnt$ and $1 \le a \le \flnt$, we have $c(u_{1+d-a}) + c(u_{1+d+a}) = c(u_{1-d-a}) + c(u_{1-d+a})$,
meaning that $\vd(u_{1+d}) = \vd(u_{1-d})$.
\hfill \qed

\vs{4mm}
Now, let $n\ge 3$ be a prime number and suppose that $c$ is a symmetric coloring of $C_n$ with respect to a vertex of $C_n$.
By Proposition \ref{pw-code}, there are at least one pair of vertices (partners) in $C_n$ that share the same code. 
This implies that $c$ is not an ID-coloring of $C_n$, and this completes the proof of Proposition \ref{prime-id-pw2}, and hence Theorem \ref{prime-id-pw} as well.

\vs{3mm}
\section{Observations for Cycles of Non-Prime Order}

The assumption that the order $n$ of the cycle $C_n$ is a prime number is essential in Theorem \ref{prime-id-pw}.
Indeed, we shall see what can be observed when the order $n$ of the cycle $C_n$ is not a prime number.
From now on, we use the notation 
$C_n = (u_0, u_1, \cdots, u_{n-1}, u_0)$,
instead of the previously used notation
$C_n = (u_1, u_2, \cdots, u_{n}, u_1)$.

\subsection{Cycles of Non-Prime Odd Order}

First, we consider cycles of non-prime odd order,
and it turns out that Theorem \ref{prime-id-pw} no longer holds.
Namely,
we obtain the following result:
\vs{-1mm}
\bthm \label{non-prime-id-pw}
Let $n$ be a non-prime odd number. Then there exists a red-white coloring of $C_n$ that is neither an ID-coloring nor a symmetric coloring with respect to any vertex of $C_n$.
\ethm 

\vs{-1mm}
Before we prove Theorem \ref{non-prime-id-pw},
we define a useful coloring of cycles with non-prime odd order.
Let $n$ be an odd number that is not prime.
Let $pq$ be a factorization of $n$, where $p$ and $q$ are neither 1 nor $n$ (note that $p$ and $q$ are odd but not necessarily prime).
Note that $p$ and $q$ are not necessarily unique.
We define a {\bf splitting-alternating coloring (SA-coloring) with $p$ splitting vertices} of $C_n = (u_0, u_1, \cdots, u_{n-1}, u_0)$ as follows.
Using the factorization $n=pq$ that is mentioned above, we color the vertices $u_{\ell q}$ $(1 \le \ell  \le p)$ (we call them {\bf splitting vertices} of $C_n$), in a way that $u_{\ell q}$ is white if $\ell$ is odd and red if $\ell$ is even.
Note that there are two ``consecutive'' splitting vertices, $u_0=u_{pq}$ and $u_q$, that have the same color, white. 
Also, the number of white splitting vertices ($\cef p2$) is greater than the number of red splitting vertices ($\flf p2$) by one.
For non-splitting vertices $u_a$ of $C_n$, where $1 \le a \le n-1$ and $a \ne \ell q$ \  $(1 \le \ell  \le p)$, 
we assign red if $a \equiv 1, 3, \cdots, q-2$ (odd) mod $q$, while we assign white if $a \equiv 2, 4, \cdots, q-1$ (even) mod $q$.
Figure \ref{C9-C15-C25} shows the SA-colorings of $C_9$, $C_{15}$ and $C_{25}$ with 3, 3, and 5 splitting vertices, respectively (all the splitting vertices are labeled, as well as the vertices $u_1$ and $u_{n-1}$).

\bfi
\scalebox{0.65}{\input{C9-C15-C25.tex}}
\capt{The SA-colorings of $C_9$, $C_{15}$ and $C_{25}$ with 3, 3, and 5 splitting vertices}
\label{C9-C15-C25}
\efi


Now, let us prove Theorem \ref{non-prime-id-pw}.\\[2mm]
{\bf Proof of Theorem \ref{non-prime-id-pw}.} \ 
Let $C_n = (u_0, u_1, \cdots, u_{n-1}, u_0)$.
Since $n$ is odd and not prime, we may assume that $n \ge 9$. Suppose that $p\ge 3$ is the smallest factor of $n$ (except for 1) and let $q = n/p$ (hence $q \ge 3$). Note that $q$ is an odd number but may not be a prime.
We consider the SA-coloring of $C_n$ with $p$ splitting vertices.

First, we show that $\vd(u_0)= \vd(u_q)=(1, 1, \cdots, 1)$, which implies that this coloring is not an ID-coloring.
Let us consider $\vd(u_0)$.
Notice that there are $\frac{p-1}2$ pairs of splitting vertices equidistant from $u_0$ (having distance $q, 2q, \cdots, \frac{p-1}2 q$ from $u_0$).
They are $u_{\ell q}$ and $u_{-\ell q}$ for each $\ell$ with $1 \le \ell \le \frac{p-1}2$.
Observe that $-\ell q = n - \ell q = pq - \ell q = (p-\ell)q$ for $1 \le \ell \le \frac{p-1}2$, and $\ell$ is odd if and only if $p-\ell$ is even.
Thus, each pair of vertices $u_{\ell q}$ and $u_{-\ell q}=u_{(p-\ell)q}$ that have distance $\ell q$ from $u_0$ contain one red vertex and one white vertex. 
Therefore, the $(\ell q)$-th coordinate of the code $\vd(u_0)$ is 1 for $1 \le \ell \le \frac{p-1}2$.
For the other coordinates of the code $\vd(u_0)$, let us consider a pair of vertices $u_b$ and $u_{-b}$ for $1 \le b \le \flnt$.
Observe that $-b \equiv q-b$ mod $q$, and hence $b \equiv m$ mod $q$ if and only if $-b \equiv q-m$ mod $q$ for ($1\le m \le q-1$).
In particular, $b \equiv m$ is odd mod $q$ if and only if $-b \equiv q-m$ is even mod $q$ for $1\le m \le q-1$.
Therefore, $c(u_b)+c(u_{-b}) = 1$ for all $1 \le b \le \flnt$ and thus
$\vd(u_0)=(1, 1, \cdots, 1)$.


Next we consider $\vd(u_q)$.
From the observation above, we obtain $q+b \equiv m$ is odd mod $q$ if and only if $q-b \equiv q-m$ is even mod $q$ for $1\le m \le q-1$, for $1 \le b \le \flnt$.
Therefore, $c(u_{q+b})+c(u_{q-b}) = 1$ for all $1 \le b \le \flnt$.
Now, there are $\frac{p-1}2$ pairs of splitting vertices equidistant from $u_q$ (having distance $q, 2p, \cdots, \frac{p-1}2 q$ from $u_q$),
which are $u_{(1+\ell)q}$ and $u_{(1-\ell)q}$ for $1 \le \ell \le \frac{p-1}2$.
When $\ell = 1$, $u_{2p}$ is red (since 2 is even) and $u_{0p}=u_n=u_{qp}$ is white (since $p$ is odd).
Observe that $(1-\ell)q = q- \ell q = n+q - \ell q = qp+q - \ell q =(p+1-\ell)q$.
Notice that $\ell$ is odd $\LRA$ $1+\ell$ is even $\LRA$ $1-\ell = p+1-\ell$ is odd.
Therefore, $c(u_{(1+\ell)q}) + c(u_{(1-\ell)q}) = 1$ for $1 \le \ell \le \frac{p-1}2$,
and hence $\vd(u_q) = (1, 1, \cdots, 1) = \vd(u_0)$.


Next, we show that this coloring is not a symmetric coloring with respect to any vertex of $C_n$, either.
What we need to show is that for any vertex $v$ of $C_n$, we can find a pair of vertices equidistant from $v$ that have different colors.
It is not difficult to see that any splitting vertex has the property, as the neighbors of them are always red and white.
So our attention goes to non-splitting vertices $u_a$ where $a\equiv 1, 2, \cdots, q-1$ mod $q$.

For the SA-colorings of $C_9$, $C_{15}$ and $C_{25}$ with 3, 3, and 5 splitting vertices respectively, we will list non-splitting vertices and corresponding pairs of vertices that are equidistant from them but have different colors. 

For $C_9= (u_0, u_1, \cdots, u_8, u_0)$, the splitting vertices are $u_0$, $u_3$ and $u_6$. For each non-splitting vertices, observe that there is a pair of vertices that are equidistant from that vertex but have different colors. Namely, 
for $u_1$, $\{u_5, u_6\}$; \ 
for $u_2$, $\{u_1, u_3\}$; \ 
for $u_4$, $\{u_2, u_6\}$; \ 
for $u_5$, $\{u_3, u_7\}$; \ 
for $u_7$, $\{u_6, u_8\}$; \ 
for $u_8$, $\{u_0, u_7\}$.
Therefore, this coloring is not symmetric.

For $C_{15}= (u_0, u_1, \cdots, u_{14}, u_0)$, the splitting vertices are $u_0$, $u_5$ and $u_{10}$. For each non-splitting vertices, observe that there is a pair of vertices that are equidistant from that vertex but have different colors. Namely, 
for $u_1$, $\{u_3, u_{14}\}$; \ 
for $u_2$, $\{u_9, u_{10}\}$; \ 
for $u_3$, $\{u_1, u_5\}$; \ 
for $u_4$, $\{u_3, u_5\}$; \ 
for $u_6$, $\{u_4, u_8\}$; \ 
for $u_7$, $\{u_4, u_{10}\}$; \ 
for $u_8$, $\{u_5, u_{11}\}$; \ 
for $u_9$, $\{u_7, u_{11}\}$; \ 
for $u_{11}$, $\{u_{10}, u_{12}\}$; \ 
for $u_{12}$, $\{u_{10}, u_{14}\}$; \ 
for $u_{13}$, $\{u_0, u_{11}\}$; \ 
for $u_{14}$, $\{u_0, u_{13}\}$.
Therefore, this coloring is not symmetric.

For $C_{25}= (u_0, u_1, \cdots, u_{24}, u_0)$, the splitting vertices are $u_0$, $u_5$, $u_{10}$, $u_{15}$ and $u_{20}$. For each non-splitting vertices, observe that there is a pair of vertices that are equidistant from that  vertex but have different colors. Namely, 
for $u_1$, $\{u_3, u_{24}\}$; \ 
for $u_2$, $\{u_9, u_{20}\}$; \ 
for $u_3$, $\{u_1, u_5\}$; \ 
for $u_4$, $\{u_3, u_5\}$; \ 
for $u_6$, $\{u_3, u_9\}$; \ 
for $u_7$, $\{u_4, u_{10}\}$; \ 
for $u_8$, $\{u_5, u_{11}\}$; \ 
for $u_9$, $\{u_7, u_{11}\}$; \ 
for $u_{11}$, $\{u_{10}, u_{12}\}$; \ 
for $u_{12}$, $\{u_{10}, u_{14}\}$; \ 
for $u_{13}$, $\{u_{11}, u_{15}\}$; \ 
for $u_{14}$, $\{u_{13}, u_{15}\}$; \ 
for $u_{16}$, $\{u_{14}, u_{18}\}$; \ 
for $u_{17}$, $\{u_{14}, u_{20}\}$; \ 
for $u_{18}$, $\{u_{15}, u_{21}\}$; \ 
for $u_{19}$, $\{u_{17}, u_{21}\}$; \ 
for $u_{21}$, $\{u_{20}, u_{22}\}$; \ 
for $u_{22}$, $\{u_{20}, u_{24}\}$; \ 
for $u_{23}$, $\{u_{0}, u_{21}\}$; \ 
for $u_{24}$, $\{u_{0}, u_{23}\}$.
Therefore, this coloring is not symmetric.

For other odd (non-prime) cycles with order $n\ge 9$, with the factorization $n=pq$ where $p$ is the smallest (but not 1) factor of $n$, we may assume that $q$ is at least 7 (again, note that $q$ is odd but may not be a prime). We consider the following cases.


(1) Suppose that $a \equiv 1, 2, \cdots, \frac{q-3}2$.
We may write that $a= \ell q + b$, where $0\le \ell \le p-1$ and $b \in \{1, 2, \cdots \frac{q-3}2\}$.
Let $d = b+1$.
Then $u_{a-d}$ is white, because $a-d = \ell q + b -(b+1) = \ell q - 1 = (\ell-1)q + q-1$.
On the other hand, $u_{a+d}$ is red, because $a+d = \ell q + b + (b + 1)= \ell q + 2b+1$ and $2b+1$ is an odd number at most $q-2$.

(2) Suppose that $a \equiv \frac{q+3}2, \frac{q+5}2, \cdots, q-1$. We may write that $a= \ell q + b$, where $0\le \ell \le p-1$ and $b \in \{ \frac{q+3}2, \frac{q+5}2,\cdots, q-1\}$.
Let $d= q+1-b$.
Now $u_{a-d}$ is white, because $a-d= \ell q +b -(q+1-b) = \ell q + 2b-1-q$, and $2b-1-q$ is an even number with $2 \le 2b-1-q \le q-3$.
On the other hand, $u_{a+d}$ is red, because $a+d = \ell q + b + (q + 1-b)= (\ell+1) q +1$.

(3) Suppose that $a \equiv \frac{q+1}2$. We may write that $a= \ell q + b$, where $0\le \ell \le p-1$ and $b = \frac{q+1}2$.
Let $d=b$.
Now, $u_{a+d}$ is red, because $a+d = \ell q + 2b = \ell q + q+1 = (\ell +1)q +1$.
On the other hand, $a-d = \ell q + b - b = \ell q$, which may or may not be white.
If $u_{\ell q}$ is white, we are done.
If $u_{\ell q}$ is red, then we consider $d' = q+b$ instead.
Now $u_{a+d'}$ is red, because $a+d' = (\ell + 2)q + 1$.
On the other hand, $u_{a-d'}$ is white, because $a-d' = (\ell -1)q$ and $u_{\ell q}$ is red.
Note that the two splitting vertices $u_0$ and $u_q$ are both white, but that does not matter in this case, although we must take caution in the following cases.
By the way, one might wonder if $d'=q+b \le \flnt$. Let us check this in detail.
Since $p\ge 3$, we have $3q \le pq = n$, i.e., $q \le \frac n3$.
Now, observe that $d' = q + b = q + \frac{q+1}2 \le \frac 32 \cdot \frac n3 + \frac 12 = \frac{n+1}2 = \cent$.
Thus, $d'$ can be greater than $\flnt$.
However, $d'$ is at most $\cent$, and fortunately, $\{ a- \flnt, a+\flnt\} = \{ a-\cent, a+\cent\}$, so there is no problem for our argument even if $d'=\cent$.

(4)
Suppose that $a = \ell q + \frac{q-1}2$ and $0 \le \ell \le p-2$. 
Let $d= \frac{q+1}2$.
Now, $u_{a-d}$ is white, because $a-d = \ell q + \frac{q-1}2 - \frac{q+1}2 = \ell q -1$.
On the other hand, for $u_{a+d}$, observe that $a+d = \ell q + \frac{q-1}2 + \frac{q+1}2 = (\ell + 1)q$, which may or may not be red.
If it is red, we are done.
If it is white, then we consider $d'=q+ \frac{q+1}2$ instead. As we saw in case (3), $d'$ is at most $\cent$, which is okay for this argument.
Now, $u_{a-d'}$ is white, since $a-d' = \ell q + \frac{q-1}2 - (q+ \frac{q+1}2) = (\ell-1)q-1$.
On the other hand, $u_{a+d'}$ is red, because $a+d' = \ell q + \frac{q-1}2 + (q+ \frac{q+1}2) = (\ell + 2)q$ and $u_{(\ell+1)q}$ is white.
Note that there are two ``consecutive'' white splitting vertices $u_0$ and $u_p$, but that case is excluded because of the assumption $\ell \ne p-1$.

(5)
Suppose that $a = (p-1)q + \frac{q-1}2$. 
Let $d = \frac{q-1}2$.
Now $u_{a-d}$ is red, since $a-d = (p-1)q$ and $p-1$ is even. 
On the other hand, $u_{a+d}$ is white, since $a+d= pq-1 = -1$.

Now the proof of the theorem is completed.  
\hfill \qed

\vs{3mm}
Combining Theorems \ref{prime-id-pw} and \ref{non-prime-id-pw}, we obtain the following statement.
\vs{-3mm}
\bthm \label{prime-id-pw-iff}
Let $n\ge3$ be an odd integer. Then the following statements are equivalent:
\vs{-4mm}
\ben
\item $n$ is a prime number;
\vs{-3mm}
\item A red-white coloring of $C_n$ is an ID-coloring if and only if it is not a symmetric coloring with respect to any vertex of $C_n$.
\een
\ethm

\subsection{Cycles of Even Order}

Next, we consider cycles of even order.
The results we have obtained so far give rise to a question as to a theorem similar to Theorem \ref{prime-id-pw} exists.
We now answer this question.

First, we need to define what it means to be symmetric for a red-white coloring of an even cycle. 
Let $n\ge 4$ be even. For a cycle $C_n$, a red-white coloring of $C_n$ is {\bf symmetric} if there is a labeling of $C_n = (u_1, u_2, \cdots, u_n, u_1)$ such that $u_k$ and $u_{n+1-k}$ have the same color for all $k$ with $1 \le k \le  \frac n2$.
Given a labeling of a symmetric coloring, each pair of vertices $u_k$ and $u_{n+1-k}$ are called {\bf partners}. Figure \ref{C10-pw} shows an example of a symmetric coloring of $C_{10}$, and each pair of partners (which have the same color) is indicated by a two-way arrow.
\vs{1mm}

\bfi
\scalebox{0.8}{
{\unitlength 0.1in%
\begin{picture}(16.7200,16.1600)(5.0500,-20.5300)%
%
\special{pn 8}%
\special{pa 1100 490}%
\special{pa 1590 493}%
\special{pa 1985 783}%
\special{pa 2134 1250}%
\special{pa 1980 1715}%
\special{pa 1583 2001}%
\special{pa 1092 1999}%
\special{pa 697 1709}%
\special{pa 548 1243}%
\special{pa 702 776}%
\special{pa 1100 490}%
\special{pa 1590 493}%
\special{fp}%
%
\special{sh 1.000}%
\special{ia 1600 501 54 54 0.0000000 6.2831853}%
\special{pn 13}%
\special{ar 1600 501 54 54 0.0000000 6.2831853}%
%
\special{sh 0}%
\special{ia 715 756 53 53 0.0000000 6.2831853}%
\special{pn 13}%
\special{ar 715 756 53 53 0.0000000 6.2831853}%
%
\special{sh 0}%
\special{ia 1110 1989 54 54 0.0000000 6.2831853}%
\special{pn 13}%
\special{ar 1110 1989 54 54 0.0000000 6.2831853}%
%
\special{sh 1.000}%
\special{ia 555 1251 50 50 0.0000000 6.2831853}%
\special{pn 13}%
\special{ar 555 1251 50 50 0.0000000 6.2831853}%
%
\special{sh 1.000}%
\special{ia 693 1703 51 51 0.0000000 6.2831853}%
\special{pn 13}%
\special{ar 693 1703 51 51 0.0000000 6.2831853}%
%
\special{sh 1.000}%
\special{ia 1984 1695 50 50 0.0000000 6.2831853}%
\special{pn 13}%
\special{ar 1984 1695 50 50 0.0000000 6.2831853}%
%
\special{sh 1.000}%
\special{ia 2126 1236 51 51 0.0000000 6.2831853}%
\special{pn 13}%
\special{ar 2126 1236 51 51 0.0000000 6.2831853}%
%
\special{sh 0}%
\special{ia 2006 780 50 50 0.0000000 6.2831853}%
\special{pn 13}%
\special{ar 2006 780 50 50 0.0000000 6.2831853}%
%
\special{sh 1.000}%
\special{ia 1110 491 54 54 0.0000000 6.2831853}%
\special{pn 13}%
\special{ar 1110 491 54 54 0.0000000 6.2831853}%
%
\special{sh 0}%
\special{ia 1560 1999 54 54 0.0000000 6.2831853}%
\special{pn 13}%
\special{ar 1560 1999 54 54 0.0000000 6.2831853}%
%
\special{pn 8}%
\special{pa 1300 770}%
\special{pa 1900 770}%
\special{fp}%
\special{sh 1}%
\special{pa 1900 770}%
\special{pa 1833 750}%
\special{pa 1847 770}%
\special{pa 1833 790}%
\special{pa 1900 770}%
\special{fp}%
%
\special{pn 8}%
\special{pa 1310 770}%
\special{pa 810 770}%
\special{fp}%
\special{sh 1}%
\special{pa 810 770}%
\special{pa 877 790}%
\special{pa 863 770}%
\special{pa 877 750}%
\special{pa 810 770}%
\special{fp}%
%
\special{pn 8}%
\special{pa 1410 1240}%
\special{pa 2010 1240}%
\special{fp}%
\special{sh 1}%
\special{pa 2010 1240}%
\special{pa 1943 1220}%
\special{pa 1957 1240}%
\special{pa 1943 1260}%
\special{pa 2010 1240}%
\special{fp}%
%
\special{pn 8}%
\special{pa 1430 1240}%
\special{pa 630 1240}%
\special{fp}%
\special{sh 1}%
\special{pa 630 1240}%
\special{pa 697 1260}%
\special{pa 683 1240}%
\special{pa 697 1220}%
\special{pa 630 1240}%
\special{fp}%
%
\special{pn 8}%
\special{pa 1310 1700}%
\special{pa 810 1700}%
\special{fp}%
\special{sh 1}%
\special{pa 810 1700}%
\special{pa 877 1720}%
\special{pa 863 1700}%
\special{pa 877 1680}%
\special{pa 810 1700}%
\special{fp}%
%
\special{pn 8}%
\special{pa 1300 1700}%
\special{pa 1900 1700}%
\special{fp}%
\special{sh 1}%
\special{pa 1900 1700}%
\special{pa 1833 1680}%
\special{pa 1847 1700}%
\special{pa 1833 1720}%
\special{pa 1900 1700}%
\special{fp}%
\end{picture}}%
}
\capt{A symmetric coloring of $C_{10}$}
\label{C10-pw}
\efi

\vs{-3mm}
Unlike symmetric colorings of odd cycles, there is no central vertex for symmetric colorings of even cycles.
Also, there is a restriction on the number of red vertices of symmetric colorings of even cycles:
a symmetric coloring of an even cycle must assign red to even number of vertices. Equivalently, if a red-white coloring assigns red to odd number of vertices of an even cycle, then the coloring is not symmetric.

Now, we present a result that is analogous to Theorem \ref{non-prime-id-pw}.

\vs{-2mm}
\bpr \label{even-not-id-pw}
Let $n \ge 4$ be an even number. Then there exists a red-white coloring of $C_n$ that is neither an ID-coloring nor a symmetric coloring with respect to any vertex of $C_n$.
\epr \vs{-3mm}
\pf \ 
Let us consider a red-white coloring $c$ that assigns red to only one vertex of $C_n$, say $v$. Since the two neighbors of $v$ in $C_n$ have the same code (the first coordinate is 1 and the other coordinates are all 0), $c$ is not an ID-coloring. Furthermore, $c$ assigns red to only one vertex, so it cannot be a symmetric coloring of $C_n$, given that there are odd number of white vertices.
\hfill \qed

\vs{2mm}
With the aid of Proposition \ref{even-not-id-pw}, we are now ready to state a generalized version of Theorem \ref{prime-id-pw-iff}.
\vs{-3mm}
\bthm \label{prime-id-pw-all}
Let $n\ge3$ be an integer. Then the following statements are equivalent:
\vs{-4mm}
\ben
\item $n$ is a prime number;
\vs{-3mm}
\item A red-white coloring $c$ of $C_n$ is an ID-coloring if and only if $c$ is not a symmetric coloring with respect to any vertex of $C_n$.
\een
\ethm


\section{Codes of Symmetric Colorings of Cycles}

We saw the following result back in Section 2.6.
\\[3mm]
{\bf Proposition \ref{pw-code} }
A red-white coloring of an odd cycle $C_n$ is a symmetric coloring with respect to some vertex of $C_n$ if and only if each pair of partners of $C_n$ have the same code.

\vs{2mm}
Given a symmetric coloring of an odd cycle $C_n$,
Proposition \ref{pw-code} answered the question as to whether a pair of partners (which share the same color) share the same code (the answer was yes).
Now we consider whether distinct pairs of partners have distinct codes. The answer is clearly yes if we consider two pairs of partners with different colors.
However, it is not straightforward to answer the question if we consider two pairs of partners with the same color.

Here we present an interesting example of a symmetric coloring of $C_{21}=(u_0, u_1, \cdots, u_{20}, u_0)$ in Figure \ref{C21-pw}.
Note that this is a symmetric coloring of $C_{21}$ with respect to the vertex $u_0$.
Since $u_3$ and $u_{18}$ are partners, they share the same color and the same code. There is another pair of partners, $u_4$ and $u_{17}$, sharing the same color and code.
Now, one may have a question as to whether $u_3$ and $u_4$ have distinct codes, given that they are not partners with respect to the vertex $u_0$.
Interestingly, the answer is no. They do share the same code $\vd(u_3)=\vd(u_4)=[1, 0, 1, 1, 0, 1, 2, 1, 0, 1]$. 
Consequently, the partners of these vertices must have the same code, so we have $\vd(u_3)=\vd(u_4)=\vd(u_{17})=\vd(u_{18})$.
What makes it happen is the multiple symmetries of this red-white coloring.
If we look carefully, we can see that this coloring is also a symmetric coloring with respect to the vertex $u_7$, and even a symmetric coloring with respect to the vertex $u_{14}$ as well!
Thus, it is undoubtedly clear that $u_3$ and $u_4$ have the same code, since they are partners with respect to a different central vertex, $u_{14}$.
Also, $u_{17}$ and $u_{18}$ have the same code, since they are partners with respect to a different central vertex, $u_7$.
Given that this coloring can be seen as a symmetric coloring with respect to three different central vertices, it turns out that there are only four distinct codes for the vertices of $C_{21}$ with respect to this coloring.
First, for the red vertices, there are only two distinct codes: $\vd(u_3)=\vd(u_4)=\vd(u_{10})=\vd(u_{11})
=\vd(u_{17})=\vd(u_{18})=[1, 0, 1, 1, 0, 1, 2, 1, 0, 1]$
and $\vd(u_0)=\vd(u_7)=\vd(u_{14})=[0,0,2,2,0,0,2,0,0,2]$.
For the white vertices, there are also only two distinct codes, $\vd(u_1) = \vd(u_6) = \vd(u_8) = \vd(u_{13}) = \vd(u_{15}) = \vd(u_{20}) = [1, 1, 1, 1, 1, 1, 0, 1, 1, 1]$
and
$\vd(u_2) = \vd(u_5) = \vd(u_9) = \vd(u_{12}) = \vd(u_{16}) = \vd(u_{19}) = [1,2,0,0,2,1,0,1,2,0]$.

\bfi
\scalebox{0.7}{
{\unitlength 0.1in%
\begin{picture}(37.8100,37.9400)(6.6000,-45.5400)%
%
\special{pn 8}%
\special{pa 2400 4510}%
\special{pa 2920 4510}%
\special{pa 3417 4357}%
\special{pa 3847 4064}%
\special{pa 4171 3657}%
\special{pa 4361 3173}%
\special{pa 4400 2655}%
\special{pa 4284 2148}%
\special{pa 4024 1697}%
\special{pa 3643 1344}%
\special{pa 3174 1118}%
\special{pa 2660 1041}%
\special{pa 2146 1118}%
\special{pa 1677 1344}%
\special{pa 1296 1697}%
\special{pa 1036 2148}%
\special{pa 920 2655}%
\special{pa 959 3173}%
\special{pa 1149 3657}%
\special{pa 1473 4064}%
\special{pa 1903 4357}%
\special{pa 2400 4510}%
\special{pa 2920 4510}%
\special{fp}%
%
\special{sh 1.000}%
\special{ia 2670 1050 54 54 0.0000000 6.2831853}%
\special{pn 13}%
\special{ar 2670 1050 54 54 0.0000000 6.2831853}%
%
\special{sh 1.000}%
\special{ia 1280 1700 54 54 0.0000000 6.2831853}%
\special{pn 13}%
\special{ar 1280 1700 54 54 0.0000000 6.2831853}%
%
\special{sh 1.000}%
\special{ia 1030 2160 54 54 0.0000000 6.2831853}%
\special{pn 13}%
\special{ar 1030 2160 54 54 0.0000000 6.2831853}%
%
\special{sh 1.000}%
\special{ia 1160 3670 54 54 0.0000000 6.2831853}%
\special{pn 13}%
\special{ar 1160 3670 54 54 0.0000000 6.2831853}%
%
\special{sh 1.000}%
\special{ia 2380 4500 54 54 0.0000000 6.2831853}%
\special{pn 13}%
\special{ar 2380 4500 54 54 0.0000000 6.2831853}%
%
\special{sh 1.000}%
\special{ia 2950 4500 54 54 0.0000000 6.2831853}%
\special{pn 13}%
\special{ar 2950 4500 54 54 0.0000000 6.2831853}%
%
\special{sh 1.000}%
\special{ia 4160 3670 54 54 0.0000000 6.2831853}%
\special{pn 13}%
\special{ar 4160 3670 54 54 0.0000000 6.2831853}%
%
\special{sh 1.000}%
\special{ia 4270 2160 54 54 0.0000000 6.2831853}%
\special{pn 13}%
\special{ar 4270 2160 54 54 0.0000000 6.2831853}%
%
\special{sh 1.000}%
\special{ia 4030 1700 54 54 0.0000000 6.2831853}%
\special{pn 13}%
\special{ar 4030 1700 54 54 0.0000000 6.2831853}%
%
\special{sh 0}%
\special{ia 2140 1120 51 51 0.0000000 6.2831853}%
\special{pn 13}%
\special{ar 2140 1120 51 51 0.0000000 6.2831853}%
%
\special{sh 0}%
\special{ia 3200 1120 51 51 0.0000000 6.2831853}%
\special{pn 13}%
\special{ar 3200 1120 51 51 0.0000000 6.2831853}%
%
\special{sh 0}%
\special{ia 3650 1360 51 51 0.0000000 6.2831853}%
\special{pn 13}%
\special{ar 3650 1360 51 51 0.0000000 6.2831853}%
%
\special{sh 0}%
\special{ia 1650 1360 51 51 0.0000000 6.2831853}%
\special{pn 13}%
\special{ar 1650 1360 51 51 0.0000000 6.2831853}%
%
\special{sh 0}%
\special{ia 930 2620 51 51 0.0000000 6.2831853}%
\special{pn 13}%
\special{ar 930 2620 51 51 0.0000000 6.2831853}%
%
\special{sh 0}%
\special{ia 4390 2620 51 51 0.0000000 6.2831853}%
\special{pn 13}%
\special{ar 4390 2620 51 51 0.0000000 6.2831853}%
%
\special{sh 0}%
\special{ia 960 3150 51 51 0.0000000 6.2831853}%
\special{pn 13}%
\special{ar 960 3150 51 51 0.0000000 6.2831853}%
%
\special{sh 0}%
\special{ia 4360 3150 51 51 0.0000000 6.2831853}%
\special{pn 13}%
\special{ar 4360 3150 51 51 0.0000000 6.2831853}%
%
\special{sh 0}%
\special{ia 3860 4040 51 51 0.0000000 6.2831853}%
\special{pn 13}%
\special{ar 3860 4040 51 51 0.0000000 6.2831853}%
%
\special{sh 0}%
\special{ia 1470 4060 51 51 0.0000000 6.2831853}%
\special{pn 13}%
\special{ar 1470 4060 51 51 0.0000000 6.2831853}%
%
\special{sh 0}%
\special{ia 1880 4330 51 51 0.0000000 6.2831853}%
\special{pn 13}%
\special{ar 1880 4330 51 51 0.0000000 6.2831853}%
%
\special{sh 0}%
\special{ia 3450 4330 51 51 0.0000000 6.2831853}%
\special{pn 13}%
\special{ar 3450 4330 51 51 0.0000000 6.2831853}%
\put(25.7000,-8.9000){\makebox(0,0)[lb]{$u_0$}}%
\put(20.9000,-10.3000){\makebox(0,0)[lb]{$u_1$}}%
\put(32.6000,-10.2000){\makebox(0,0)[lb]{$u_{20}$}}%
\put(10.3000,-37.4000){\makebox(0,0)[rt]{$u_7$}}%
\put(42.5000,-37.6000){\makebox(0,0)[lt]{$u_{14}$}}%
\put(11.8000,-16.6000){\makebox(0,0)[rb]{$u_3$}}%
\put(9.6000,-21.1000){\makebox(0,0)[rb]{$u_4$}}%
\put(41.4000,-16.5000){\makebox(0,0)[lb]{$u_{18}$}}%
\put(43.4000,-21.0000){\makebox(0,0)[lb]{$u_{17}$}}%
\end{picture}}%
}
\capt{A symmetric coloring of $C_{21}$}
\label{C21-pw}
\efi

\vs{1mm}
This example gives rise to the following result on symmetric colorings of cycles with non-prime odd order.

\vs{-2mm}

\bthm \label{nonprime-pw-not-k-pair}
Let $n\ge 9$ be an odd integer that is not prime. 
Then there exists a red-white coloring of the cycle $C_n$
that is a symmetric coloring with respect to more than one vertex.
\ethm

In order to consider Theorem \ref{nonprime-pw-not-k-pair}, it is useful to observe the following, which is obtained directly from the definition of symmetric colorings of cycles.

\bobs \label{code02}
Let $n\ge 3$ be an odd integer and let $c$ be a symmetric coloring of the cycle $C_n$. Then the vertex $u \in V(C_n)$ is a central vertex of $c$ if and only if the code $\vd(u)$ does not contain 1.
\eobs
\vs{-1mm}
{\bf Proof of Theorem \ref{nonprime-pw-not-k-pair}. }\hs{-1mm}
Since $n$ is not prime, let $n=pq$ where $p$ and $q$ are positive odd integers that are neither 1 nor $n$ (but not necessarily prime). 
Note that $p, q \ge 3$.
Let $c$ be a red-white coloring of the cycle $C_n = (u_0, u_1, u_2, \cdots, u_{n-1}, u_0)$ where the vertices $u_{\ell p}$ are red for $0 \le \ell \le q-1$ and the rest of the vertices are white. 
For each $0 \le \ell \le q-1$, $\vd(u_{\ell p})$ is the code with $\flnt$ entries, where the $(kp)$-th coordinate is 2 for $1 \le k \le \frac{q-1}2$ and the other coordinates are all 0.
Thus, this coloring is a symmetric coloring with respect to exactly $q$ vertices by Observation \ref{code02}. 
\hfill \qed

\vs{3mm}

In spite of what we have just obtained, we can observe a totally different result when the order of the cycle $C_n$ is prime.
First, we prove the following.

\vs{-3mm}
\bthm \label{prime-pw-central}
Let $n \ge 3$ be a prime number. 
Then any symmetric coloring of the cycle $C_n$ that is neither all-white nor all-red has exactly one central vertex.
\ethm
\vs{-3mm}
\pf \ 
Let $c$ be a symmetric coloring of the cycle $C_n = (u_0, u_1, \cdots, u_{n-1}, u_0)$.
We may assume that $u_0$ is a central vertex of the coloring $c$.
Assume to the contrary that there exists another central vertex $u_d$ of $c$, where $1 \le d \le \flnt$.
First, suppose that $u_0$ is white.
If $u_d$ is red, then the sum of all the coordinates of the code $\vd(u_d)$ must be odd, since the code refers to odd number of red vertices. Thus, $\vd(u_d)$ must contain 1 in some coordinate, which contradicts Observation \ref{code02}.
Now we suppose that $u_d$ is white.
Note that $u_{2d}$ is the partner of $u_0$ with respect to the vertex $u_d$, so $u_0$ and $u_{2d}$ share the same color and code. Namely, $u_{2d}$ is white and $\vd(u_{2d}) = \vd(u_0)$. Since $\vd(u_{2d})$ does not contain 1, $u_{2d}$ is another central vertex of $C_n$. Next, it follows that $u_{3d}$ and $u_d$ are partners with respect to the vertex $u_{2d}$, and hence $u_{3d}$ is white and $\vd(u_{3d})=\vd(u_d)$, which contains only 0 and/or 2, making $u_{3d}$ another central vertex of $C_n$. Inductively, it follow that $u_{kd}$ is white for any positive integer $k$. 
By algebra, $\{ kd \in \ZZ/n\ZZ\  |\  k \in \NN \} = \ZZ/n\ZZ$ and hence $\{ u_{kd} \in V(C_n) \ | \ k \in \NN \} = V(C_n)$, given that $n$ is a prime number.
Therefore, all the vertices of $C_n$ are colored white, which contradicts the assumption.
The remaining case ($u_0$ is red) is shown using exactly the same argument.
\hfill \qed

\vs{3mm}

Theorem \ref{prime-pw-central} leads us to the converse of Theorem \ref{nonprime-pw-not-k-pair}.

\vs{-1mm}
\bthm \label{prime-pw-flnt}
Let $n \ge 3$ be a prime number.
For any symmetric coloring of the cycle $C_n$ that is neither all-white nor all-red, distinct pairs of partners have distinct codes. 
\ethm
\vs{-3mm}
\pf \ 
Let $c$ be a symmetric coloring of the cycle $C_n$ that is neither all-white nor all-red, with the central vertex $v$.
Assume to the contrary that there exist two distinct vertices $x$ and $y$ that are not partners with respect to $v$ such that $\vd(x) = \vd(y)$.
Using Algorithm \ref{thealgorithm}, we obtain a symmetric coloring $c'$ of $C_n$, which must coincide the given coloring $c$.
Let $v'$ be the central vertex of $c'$.
Since $x$ and $y$ are partners with respect to $v'$ but they are not with respect to $v$,
it follows that $v\ne v'$. 
This contradicts Theorem \ref{prime-pw-central}. 
\hfill \qed

\vs{3mm}

Combining Theorems \ref{nonprime-pw-not-k-pair}, \ref{prime-pw-central} and \ref{prime-pw-flnt}, we obtain the following.
\vs{-2mm}
\bthm 
For a positive odd integer $n\ge 3$, the following statements are equivalent:
\\
(a) $n$ is a prime number;
\\
(b) Any symmetric coloring of the cycle $C_n$ that is neither all-white nor all-red has exactly one central vertex;
\\
(c) For any symmetric coloring of the cycle $C_n$ that is neither all-white nor all-red, distinct pairs of partners have distinct codes. 
\ethm


\section{Remarks}

\vs{-2mm}

In this paper, we established a criterion to determine whether a red-white coloring of a cycle with a prime order is an ID-coloring or not.
Having such a criterion is useful, not only because it can be used for cycles, but also because it can be applied for graphs that contain a cycle as a subgraph, as explained in the introduction.
The same thing can be said for a graph $H$ in general: if we can establish a criterion to determine whether a red-white coloring of a graph $H$ is an ID-coloring or not, then it will be useful, not only because it can be used for $H$ itself, but also because it can be applied for graphs that contain $H$ as a subgraph.
Therefore, one of our next goals is to find more criteria for various classes of graphs to determine whether red-white colorings of them are ID-coloring or not.

\vs{5mm}
\medskip

\noindent
{\bf Acknowledgment:}  We are grateful to Professor Gary Chartrand for suggesting concepts and problems to us and kindly providing useful information on this topic. 
Also we are thankful to Professor Ping Zhang for reading the early draft of Section 2 (that we obtained in 2022) and giving us a positive feedback.

\vs{2mm}
Author Information:\\[3mm]
Yuya Kono 
\\[3mm]
Gakushuin Boys' Senior High School\\[1mm]
1-5-1 Mejiro, Toshima-ku, Tokyo\\
171-0031, Japan
\\[2mm]
20240235@gakushuin.ac.jp\\[1mm]
yuya.kono@wmich.edu

\end{document}